\documentclass{notices}
\usepackage{amsfonts,amssymb,amsmath,amscd}
\usepackage{graphicx,xcolor,comment}\def\d{{\, \rm d}}

\usepackage{listings}
\usepackage{color} %red, green, blue, yellow, cyan, magenta, black, white
\definecolor{mygreen}{RGB}{28,172,0} % color values Red, Green, Blue
\definecolor{mylilas}{RGB}{170,55,241}

\title{
Taming Uncertainty in a Complex World: The Rise of Uncertainty Quantification --- A Tutorial for Beginners
}

\author{
  Nan Chen
  \affil{
    Nan Chen is an Associate Professor of Mathematics at the University of Wisconsin-Madison. His email address is chennan@math.wisc.edu.
    }
  \and
  Stephen Wiggins
  \affil{
    Stephen Wiggins is the William R. Davis ’68 Distinguished Chair in Mathematics at the United State Naval Academy and Professor of Applied Mathematics at the University of Bristol. His email address is s.wiggins@bristol.ac.uk.
   }
   \and
   Marios Andreou
  \affil{
    Marios Andreou is a PhD student of Mathematics at the University of Wisconsin-Madison. His email address is mandreou@math.wisc.edu.
   }
}
\date{}

\begin{document}

\maketitle

\section*{}
George Box, a British statistician, wrote the famous aphorism, ``All models are wrong, but some are useful.'' The aphorism acknowledges that models, regardless of qualitative, quantitative, dynamical, or statistical, always fall short of the complexities of reality. The ubiquitous imperfections of models come from various sources, including the lack of a perfect understanding of nature, limited spatiotemporal model resolutions due to computational power, inaccuracy in the initial and boundary conditions, etc. Therefore, uncertainty quantification (UQ), which quantitatively characterizes and estimates uncertainties, is essential to identify the usefulness of the model. UQ also defines the range of possible outcomes when certain aspects of the system are not precisely known. A key objective of UQ is to explore how uncertainties propagate, both through time evolution and across different quantities via complex nonlinear dependencies.

\section*{Characterizing Uncertainties}
When uncertainty appears in the model and the input, the output can potentially take different values with typically unequal chances. Therefore, it is natural to characterize the output as a random variable, and the UQ of the model output can be based on the associated probability density function (PDF).

Intuitively, the spread of a distribution, which describes how close the possible values of the output are to each other, measures the uncertainty of a variable. If a PDF is Gaussian, variance is a natural indicator describing the spread. Therefore, the uncertainty associated with the three Gaussian distributions is expected to increase from Panel (a) to Panel (c) in Figure \ref{DifferentPDFs}. On the other hand, non-Gaussian distributions are widely seen in practice due to the intrinsic nonlinearity of the underlying system. The non-Gaussian features should also be highlighted in the rigorous quantification of the uncertainty.

\subsection*{Shannon's entropy: quantifying the uncertainty in one PDF}
Denote by $p(x)$ the PDF of a random variable $X$. Shannon's entropy, an information measurement, is a natural choice to rigorously quantify uncertainty. It is defined as \cite{gray2011entropy}:
\begin{equation}\label{Shannon_Entropy}
    \mathcal{S}(p) = -\int p(x) \ln(p(x)) \d x.
  \end{equation}
Shannon's entropy originates from the theory of communication when a ``word'' is represented as a sequence of binary digits with length $n$, so the set of all words of length $n$ has $2^n=N$ elements. Therefore, the amount of information needed to characterize one element, which is the number of digits, is $n=\log_2N$. Consider the case where the entire set is divided into disjoint subsets, each with $N_i$ total elements. The chance of randomly taking one element that belongs to the $i$-th subset is $p_i=N_i/N$. If an element belongs to the $i$-th subset, then the additional information to determine it is $\log_2N_i$. Therefore, the average amount of information to determine an element is
\begin{equation}
\begin{split}
    \sum_i { \frac{N_i}{N}}{ \log_2 N_i} &=  \sum_i { \frac{N_i}{N}}{ \log_2\left( {\frac{N_i}{N}}\cdot N\right)}\\&= \sum_i p_i\log_2p_i + \log_2N.
\end{split}
\end{equation}
Recall that $\log_2N$ is the information to determine an element given the full set. Thus, the corresponding average lack of information is $-\sum p_i \log_2p_i$, which is the uncertainty. The formal definition of Shannon's entropy \eqref{Shannon_Entropy} generalizes the above argument. It exploits $-\ln(p(x))$, the negative of a natural logarithm function, to characterize the lack of information of each event $x$ and then takes the continuous limit to replace the finite summation with an integral that represents the uncertainty averaged over all events. See \cite{shannon1948mathematical} for a more rigorous derivation of Shannon's entropy and its uniqueness. Notably, Shannon's entropy \eqref{Shannon_Entropy} applies to general non-Gaussian PDFs.

\begin{figure}[htbp]
	\begin{center}
		\hspace*{-0cm}\includegraphics[width=8cm]{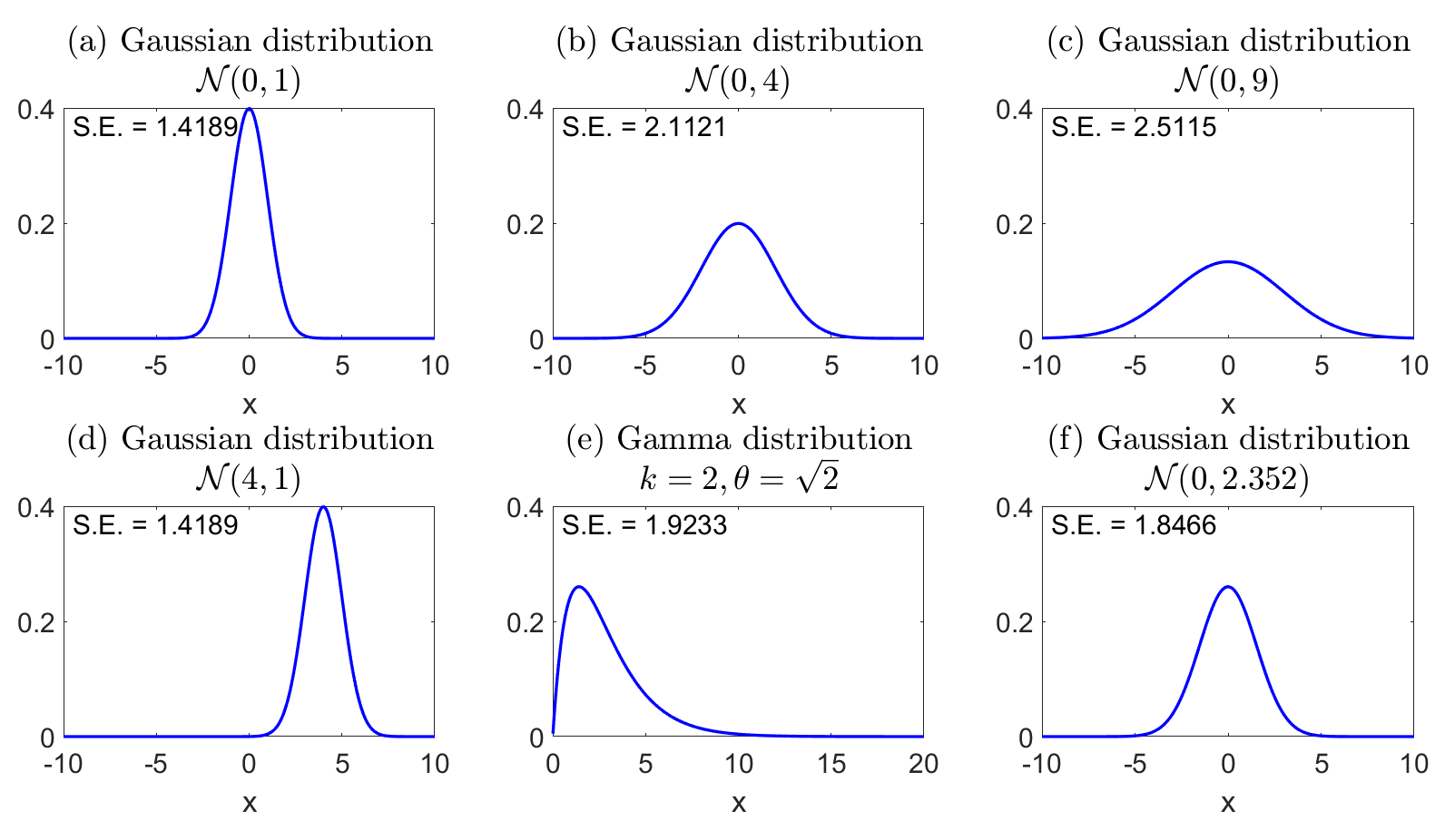}
	\end{center}
	\caption{Quantifying the uncertainty using Shannon's entropy (S.E.) for Gaussian and non-Gaussian distributions. In all panels, the x-axis spans $20$ units.}
	\label{DifferentPDFs}
\end{figure}

For certain distributions, $\mathcal{S}(p)$ can be written down explicitly. If $p\sim\mathcal{N}(\mu, R)$ is a $m$-dimension Gaussian distribution, where $\mu$ and $R$ are the mean and covariance, then Shannon entropy has the following form:
  \begin{equation}\label{Shannon_Entropy_Gaussian}
    \mathcal{S}(p) = \frac{m}{2}(1+\ln2\pi) + \frac{1}{2}\ln\det(R).
  \end{equation}
In the one-dimensional situation, \eqref{Shannon_Entropy_Gaussian} implies that the uncertainty is uniquely determined by the covariance and is independent of the mean value, consistent with the intuition. Therefore, Shannon's entropy confirms that the uncertainty increases from Panel (a) to Panel (c) in Figure \ref{DifferentPDFs}, and the uncertainty has the same value in Panel (a) and (d). However, without a systematic information-based measurement, there is typically no single empirical indicator, like variance (or equivalently, peak height), that can fully characterize the uncertainty in complex non-Gaussian PDFs. Therefore, Shannon's entropy provides a systematic way for UQ. Panel (e) shows a Gamma distribution, which is skewed and has a one-sided fat tail. The fat tail usually corresponds to extreme events, which are farther from the mean value, naturally increasing the uncertainty. Therefore, despite the same peak height of the two PDFs in Panels (e) and (f), the fat-tailed Gamma PDF has a larger entropy.

\subsection*{Relative entropy: measuring how one PDF is different from another}

In many practical problems, the interest lies in estimating the lack of information in one distribution $p^M(x)$ related to another $p(x)$. Typically, $p(x)$ is associated with a full probabilistic model while $p^M(x)$ comes from a reduced-order approximate model. The latter is often less informative but is widely used in practice to accelerate computations. %Another situation is that $p^M(x)$ represents the uncertainty associated with a given model while $p(x)$ is obtained by combining the information from the same model with the additional observational data. The goal is to quantify the uncertainty reduction in $p(x)$ with the help of observations.

Recall that the uncertainty about $x$ in $p(x)$ is $\mathcal{I}(p(x)) = -\ln(p(x))$. The average lack of information, i.e., Shannon's entropy, is
      \begin{equation}\label{RE_Entropy_PM}
        \langle \mathcal{I}(p(x)) \rangle = -\int p(x) \ln(p(x)) \d x.
      \end{equation}
However, while the uncertainty about $x$ in $p^M(x)$ is $\mathcal{I}(p^M(x)) = -\ln (p^M(x))$, the expected lack of information under $p^M(x)$ is
   \begin{equation}\label{RE_Entropy_IM}
        \langle \mathcal{I}(p^M(x)) \rangle = -\int {p(x)} \ln(p^M(x)) \d x.
   \end{equation}
This is because even though the approximate model is used to measure the information $\ln(p^M(x))$, the actual probability of $x$ to appear always comes from the full system, namely $p(x)$, which is objective and independent of the choice of the model. In other words, the role of the approximate model is to provide the lack of information for each event $x$. In contrast, the underlying distribution of the occurrence of $x$ is objective regardless of the approximate model used. Relative entropy, also known as the Kullback-Leibler (KL) divergence, characterizes the difference between these two entropies \cite{cover1999elements}:
  \begin{equation}\label{Relative_Entropy}
  \begin{split}
    \mathcal{P}(p,p^M)&=\langle \mathcal{I}(p^M(x)) \rangle - \langle \mathcal{I}(p(x)) \rangle\\
     &= \int p(x) \ln\left(\frac{p(x)}{p^M(x)}\right) \d x.
  \end{split}
  \end{equation}
The relative entropy $\mathcal{P}(p,p^M)$ is non-negative. It becomes larger when $p$ and $p^M$ are more distinct. Since $p(x)$ appears in both $\langle \mathcal{I}(p^M(x)) \rangle$ and $\langle \mathcal{I}(p(x)) \rangle$, the relative entropy is not symmetric. One desirable feature of $\mathcal{P}(p,p^M)$ is that it is invariant under general nonlinear change of variables. As a remark, if ${p^M(x)} \ln(p^M(x))$ is utilized in \eqref{RE_Entropy_IM} to compute the expected lack of information under $p^M(x)$, then the analog to \eqref{Relative_Entropy} is called Shannon entropy difference.

To illustrate that relative entropy \eqref{Relative_Entropy} is a more appropriate definition of the lack of information than Shannon entropy difference, consider a model to bet on a soccer game: Team A vs. Team B. A comprehensive model $p$ gives the odds for Team A: 10\% (win), 10\% (draw), and 80\% (lose). However, someone unfamiliar with soccer may have a biased model $p^M$, which gives the odds for Team A: 80\% (win), 10\% (draw), and 10\% (lose). If entropy difference is used, the resulting lack of information will be precisely zero. In contrast, relative entropy considers the lack of information in $p^M$ related to $p$ regarding each outcome (win/draw/lose).

When both $p\sim\mathcal{N}(\mu,R)$ and $p^M\sim\mathcal{N}(\mu^M,R^M)$ are $m$-dimensional Gaussians, the relative entropy has the following explicit formula \cite{majda2002mathematical},
\begin{equation}\label{Relative_Entropy_Gaussian}
\begin{split}
    \mathcal{P}&(p,q)= \textstyle\frac{1}{2}\Big[(\mu-\mu^M)^\mathtt{T}(R^M)^{-1}(\mu-\mu^M)\Big] \\&+ \textstyle\frac{1}{2}\Big[\mbox{tr}(R(R^M)^{-1})-m-\ln\mbox{det}(R(R^M)^{-1})\Big],
\end{split}
\end{equation}
where 'tr' and 'det' are the trace and determinant of a matrix, respectively. The first term on the right-hand side of \eqref{Relative_Entropy_Gaussian} is called `signal', which measures the lack of information in the mean weighted by model covariance. The second term involving the covariance ratio is called `dispersion'.

\section*{UQ in Dynamical Systems}
From now on, UQ will be discussed in the context of dynamical systems. The uncertainty propagates in different ways in linear and nonlinear systems.
\subsection*{Examples of uncertainty propagation in linear and nonlinear dynamical systems}
\begin{figure}[htbp]
	\begin{center}
		\hspace*{-0cm}\includegraphics[width=8cm]{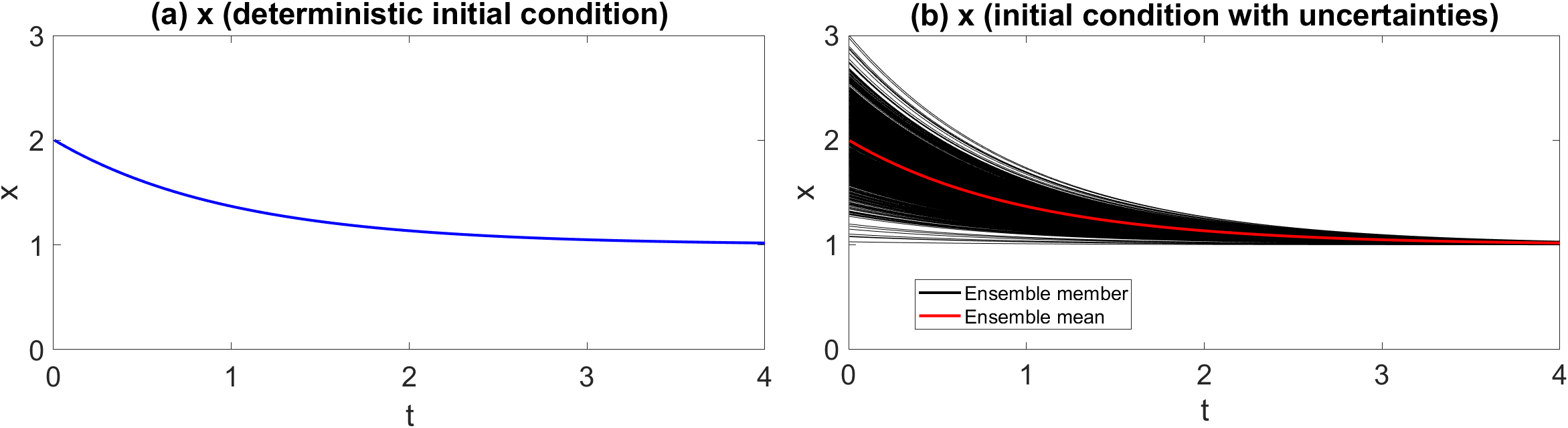}
	\end{center}
	\caption{Solutions of the linear system \eqref{Linear_ODE}. Panel (a): time evolution of $x(t)$ with a deterministic initial condition. Panel (b): time evolution of $x(t)$ with the initial condition given by a Gaussian distribution. The ensemble size is 1000. }
	\label{ODE_UQ}
\end{figure}

Consider a linear ordinary differential equation,
\begin{equation}\label{Linear_ODE}
\frac{\d x}{\d t} = -a x + f,
\end{equation}
where $a>0$ is the damping coefficient and $f$ is an external forcing. The solution of \eqref{Linear_ODE} can be written down explicitly $x(t) = x(0)e^{-at} +(1-e^{-at})f/a$, where $x(t)$ converges to $f/a$ in an exponential rate. Consider the time evolution of the linear dynamics \eqref{Linear_ODE} with two different initial conditions. A deterministic initial condition $x(0) = 2$ is given in the first case. The time evolution of $x(t)$ is shown in Panel (a) of Figure \ref{ODE_UQ}. In the second case, uncertainty appears in the initial condition, which is given by a Gaussian distribution $x(0) \sim \mathcal{N}(2,0.09)$. Different initial values are drawn from this distribution and follow the governing equation \eqref{Linear_ODE}. The black curves in Panel (b) show the time evolution of different ensemble members, while the red curve is the ensemble average. The uncertainty dissipated over time, and the time evolution of the ensemble average follows the same trajectory as the deterministic case. In other words, the uncertainty does not change the mean dynamics.

\begin{figure}[htbp]
	\begin{center}
		\hspace*{-0cm}\includegraphics[width=8cm]{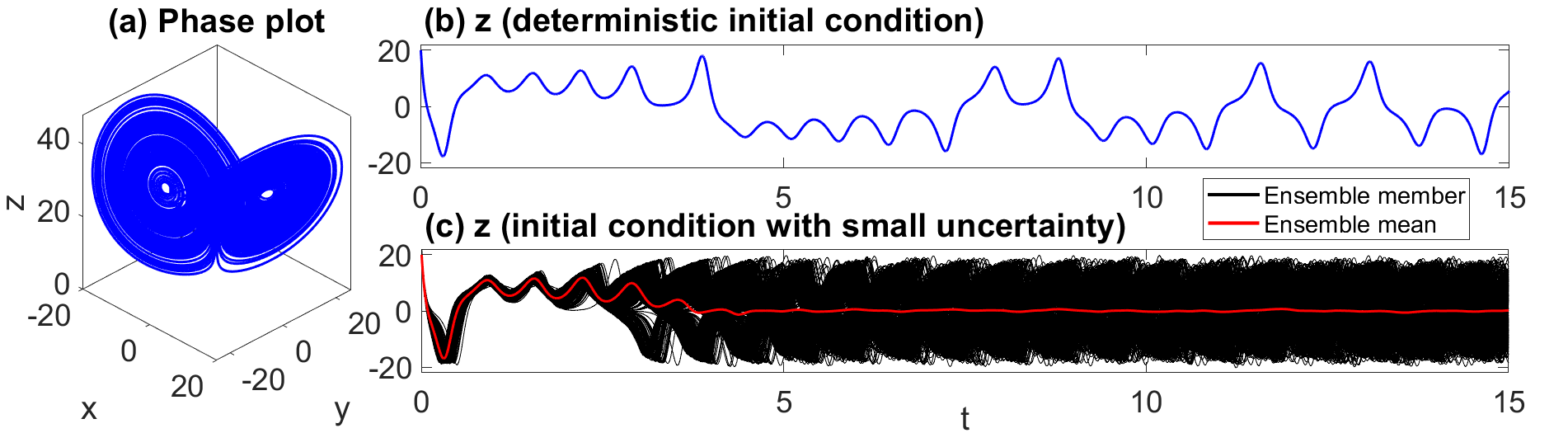}
	\end{center}
	\caption{Solutions of the nonlinear chaotic Lorenz 63 model \eqref{Lorenz63}. Panel (a): the Lorenz attractor. Panel (b): time evolution of $z(t)$ with a deterministic initial condition. Panel (b): time evolution of $z(t)$ with the initial condition given by a Gaussian distribution. The ensemble size is 1000.  }
	\label{L63_UQ}
\end{figure}

Next, consider the Lorenz 63 model, which is a nonlinear chaotic system \cite{lorenz1963deterministic},
\begin{equation}\label{Lorenz63}
                \frac{\d{x}}{\d t} = \sigma({y}-{x}), \quad
                \frac{\d{y}}{\d t} = {x}(\rho-{z}) - {y}, \quad
                \frac{\d{z}}{\d t} = {x}{y}-\beta {z}.
\end{equation}
Again, consider the time evolution of the solution starting from two sets of initial conditions, one deterministic $x(0) = 20, y(0)=-20, z(0) = 25$ and one containing uncertainty $x(0) \sim \mathcal{N}(20,1), y(0)\sim \mathcal{N}(-20,1), z(0) \sim \mathcal{N}(25,1)$. By taking the standard parameter values $\sigma = 10$, $\beta = 8/3$ and $\rho = 28$, the system displays a chaotic behavior. Panel (a) of Figure \ref{L63_UQ} shows the attractor of the system, which resembles a butterfly. Panel (b) shows the time evolution of $z(t)$ starting from the deterministic initial condition. In contrast, the black curves in Panel (c) show ensemble members starting from different values drawn from the given initial distribution, while the red curve is the ensemble average. Although the ensemble mean follows the trajectory in the deterministic case and has a small uncertainty within the first few units, the ensemble members diverge quickly. Notably, the ensemble average significantly differs from any model trajectories. This implies that uncertainty has a large impact on the mean dynamics. In other words, the mean evolution and the uncertainty cannot be considered separately, as in the linear case.

\subsection*{Impact of nonlinearity on uncertainty propagation}
To understand the impact of nonlinearity on uncertainty propagation, let us decompose a random variable $x$ into its mean and fluctuation parts via the Reynolds decomposition \cite{muller2006equations}
\begin{equation}\label{Reynolds_decomposition}
    x = \langle x\rangle + x^\prime,\qquad \mbox{with~} \langle x^\prime\rangle=0,
\end{equation}
where $\langle\cdot\rangle$ represents the ensemble average computed by taking the summation of all ensemble members and then dividing by the total number of ensembles. So, $\langle (x')^2\rangle$ is the variance of $x$. For the linear system \eqref{ODE_UQ}, taking the ensemble average leads to the mean dynamics
\begin{equation}
\frac{\d\langle x\rangle}{\d t} = (-a\langle x\rangle+f),
\end{equation}
which implies the time evolution of the mean $\langle x\rangle$ is not affected by the uncertainty described by the fluctuation part. However, the situation becomes very different if nonlinear terms appear in the dynamics. Consider a quadratic nonlinear term on the right-hand side of the dynamics,
\begin{equation}\label{nonlinear_ODE}
\frac{\d x}{\d t} = bx^2+ ...
\end{equation}
Since
\begin{align}
\langle x^2\rangle &= \langle(\langle x\rangle+ x^\prime)^2\rangle = \langle(\langle x\rangle^2+ (x^\prime)^2 + 2\langle x\rangle x^\prime)\rangle\notag\\
& = \langle x\rangle^2+ \langle(x^\prime)^2\rangle,\label{mean_fluctuation_variable}
\end{align}
taking the ensemble average of \eqref{nonlinear_ODE} leads to
\begin{equation}\label{nonlinear_ODE_mean_fluctuation}
  \frac{\d\langle x\rangle}{\d t} =  b\langle x\rangle^2 + b\langle (x')^2\rangle + ...,
\end{equation}
which indicates that the higher-order moments containing the information of uncertainties affect the time evolution of the lower-order moments (e.g., the mean dynamics) via nonlinearity. It also reveals that the moment equations for general nonlinear dynamics will never be closed. In practice, approximations are made to handle the terms involving higher-order moments in the governing equations of the lower-order moments to form a solvable closed system \cites{majda2018strategies, chen2023stochastic}.

In chaotic systems, small uncertainties are quickly amplified by positive Lyapunov exponents, making an accurate state forecast/estimation challenging. Additional resources can be combined with models to facilitate uncertainty reduction in state estimation.

\section*{Uncertainty Reduction via Data Assimilation (DA)}
Model and observational data are widely utilized to solve practical problems. However, neither model nor observation is close to perfect in most applications. Models are typically chaotic and involve large uncertainties. Observations contain noise and are often sparse, incomplete, and indirect. Nevertheless, when a numerical model and observations are optimally combined, the estimation of the state can be significantly improved. This is known as data assimilation (DA).
DA was initially developed in numerical weather prediction, which improves the initialization for a numerical forecast model. Since then, DA has become an essential tool for many applications, including dynamical interpolation of missing data, inferring the unobserved variables, parameter estimation, assisting control and machine learning, etc.

Denote by ${u}$ the state variable and ${v}=gu+\epsilon$ the noisy observation with $g$ the observational operator and $\epsilon$ a random noise, the underlying principle of DA is given by the Bayes' theorem \cite{jeffreys1973scientific}:
\begin{equation}\label{Bayes_formula}
{ \underbrace{p({u}|{v})}_{\mbox{posterior}}} \propto~ { \underbrace{p({u})}_{\mbox{prior}}}~\times\underbrace{p({v}|{u})}_{\mbox{likelihood}},
\end{equation}
where `$\propto$' stands for `proportional to'. In \eqref{Bayes_formula}, $p({u})$ is the forecast distribution by using a model built upon prior knowledge, while $p({v}|{u})$ is the probability of observation under the model assumption. Their combination is the conditional distribution of ${u}$ given $v$, which is called the posterior distribution. Notably, due to the additional information from observation, uncertainty is expected to be reduced from the prior to the posterior distribution.

\begin{figure}[htbp]
	\begin{center}
		\hspace*{-0cm}\includegraphics[width=8cm]{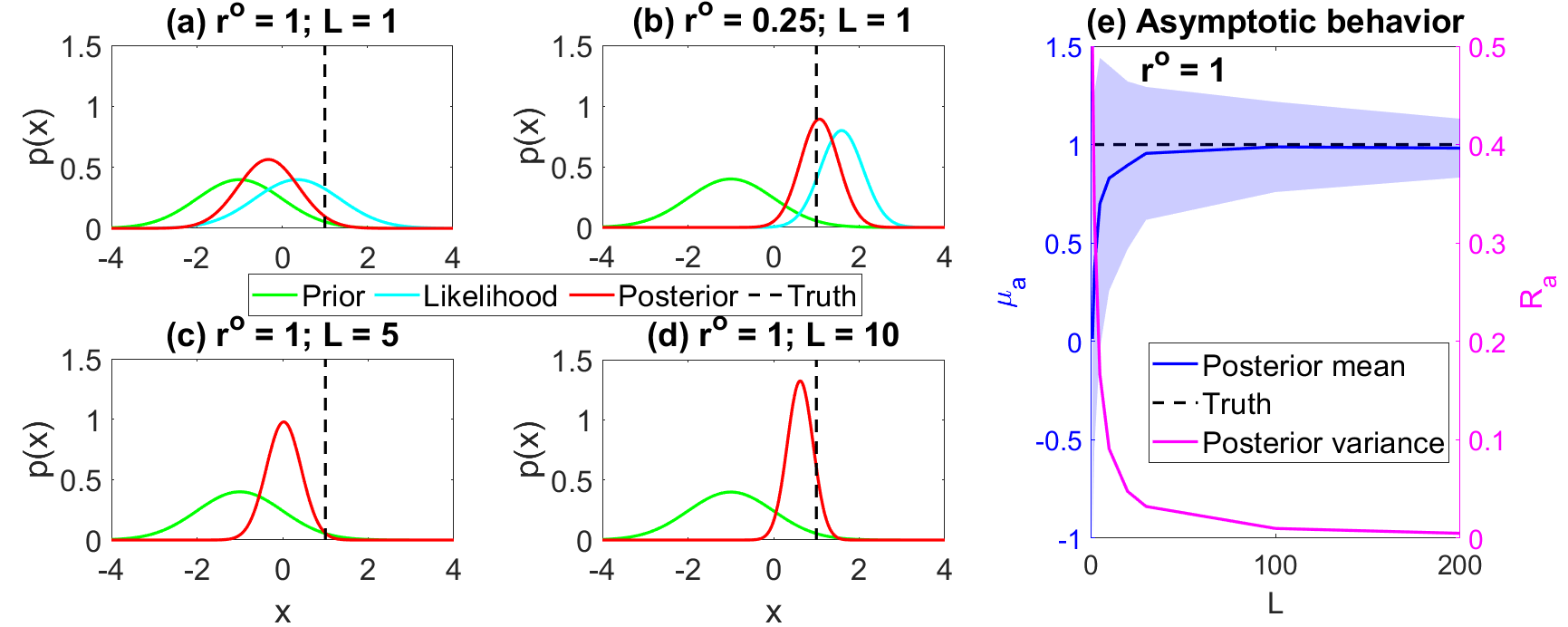}
	\end{center}
	\caption{Posterior distribution from Bayes formula \eqref{Bayes_formula} with $\mu_f=1$, $R_f=1$ and $m=1$. Panels (a)--(d): the PDFs with different observational uncertainty $r^o$ and number of observations $L$. To keep the figure concise, the PDFs corresponding to the $L=5$ or $10$ noisy observations are omitted from Panels (c)--(d). Panel (e): the asymptotic behavior of the posterior mean $\mu_a$ and the posterior variance $R_a$ as a function of $L$. Due to the randomness in observations, the shading area shows the variation of the results with 100 sets of independent observations for each fixed $L$.  }
	\label{Bayes}
\end{figure}

Assume the prior distribution $p(u)\sim\mathcal{N}(\mu_f,R_f)$ and the observational noise with variance $r^o$ are both Gaussian. Plugging these into \eqref{Bayes_formula}, it is straightforward to show that the posterior distribution $p({u}|{v})\sim\mathcal{N}(\mu_a,R_a)$ is also Gaussian. The posterior mean $\mu_a$ and covariance $R_a$ can be written down explicitly
\begin{equation}\label{Kalman_MultiD_Posterior}
%\begin{split}
  \mu_{a} = (I-Kg)\mu_{f} + K v,\quad
  R_{a} = (I - Kg) R_{f},
%\end{split}
\end{equation}
where $I$ is an identity matrix of size $m\times m$ with $m$ being the dimension of $u$ and $K$ is given by $K= R_{f}g^\mathtt{T}(gR_{f}g^\mathtt{T} + {r}^o)^{-1}$.
Now consider the case with $m=1$ and $g=1$. Then, all three distributions in \eqref{Bayes_formula} are one-dimensional Gaussians, and the observation equals the truth plus noise. In such a case, $K=R_f/(R_f+r^o)\in[0,1]$ and $\mu_{a}$ becomes a weighted summation of the prior mean and the observation with weights being $1-K$ and $K$, respectively. When the observational noise is much smaller than the prior uncertainty, i.e., $r^o \ll R_f$, the posterior mean almost fully trusts the observation since $K\approx 1$. In contrast, if the observation is highly polluted, i.e., $r^o \gg R_f$, then $K\approx 0$ and the observational information can almost be ignored, and the posterior mean nearly equals the prior mean. Essentially, the weights are fully determined by the uncertainties. On the other hand, the posterior variance $R_{a}$ in \eqref{Kalman_MultiD_Posterior} is always no bigger than the prior variance $R_f$, indicating that the observation helps reduce the uncertainty. Panels (a)--(b) of Figure \ref{Bayes} validate the above conclusions. In some applications, repeated measurements with independent noises are available, which further advance the reduction of the posterior uncertainty. For example, in the above scalar state variable case, when $L$ repeated observations are used, the observational operator becomes a $L\times 1$ vector $g=(1,\ldots,1)^\mathtt{T}$ and similarly for the noise $\epsilon = (\epsilon_1,\ldots,\epsilon_L)^\mathtt{T}$. According to Panel (a), (c), and (d), the posterior uncertainty reduces when the number of observations increases from $L=1$ to $10$. Panel (e) shows that when $L$ becomes large, the posterior mean converges to the truth, and the posterior variance decreases to zero. Correspondingly, when characterizing the uncertainty reduction in the posterior distribution related to the prior via the relative entropy \eqref{Relative_Entropy_Gaussian}, the signal part converges to a constant while the dispersion part scales as $\ln(L)$. See the online supplementary document for more details.

When the forecast model and the observational operator are linear, and the noises are Gaussian, the above procedure is called the Kalman filter \cite{kalman1960new}. Closed analytic formulae are available for the posterior distribution. However, models are generally nonlinear, and the associated PDFs are often non-Gaussian. Numerical methods, such as ensemble Kalman filter (EnKF) or particle filter, are widely used to find approximate solutions \cites{law2015data, asch2016data}.

In the following, UQ will be discussed in the context of Lagrangian DA (LaDA) \cite{apte2008bayesian}, which exploits Lagrangian observations. The same procedure for quantifying and reducing the uncertainty can be applied using other observations.

\subsection*{Lagrangian data assimilation (LaDA)}
Lagrangian tracers are moving drifters, such as robotic instruments, balloons, sea ice floes, and even litter. They are often used to recover the flow field that drives their motions. However, recovering the entire flow field based solely on Lagrangian tracers is challenging. This is because tracers are usually sparse, which prevents a direct estimation of the flow velocity in regions with no observations. In addition, the measured quantity is the tracer displacements. Observational noise can propagate through the time derivative from displacement to velocity. A flow model, despite being generally turbulent, can provide prior knowledge about the possible range of velocity in the entire domain. Observational information dominates the state estimation at the locations covered by tracers. They also serve as the constraints, conditioned on which the uncertainty in estimating the flow field through the model in regions without observations can be significantly reduced. LaDA is widely used to provide a more accurate recovery of the flow field, facilitating the study of flow properties. For simplicity, assume the tracer velocities are equal to the underlying flow field. The LaDA scheme consists of two sets of equations: one for the observational process and the other for the flow model, described as follows:
\begin{align}
  \mbox{Obs process:}\quad\frac{\d x_l}{\d t} &= u(x_l,t) + \sigma_x\dot{W}_l,\quad l =1, \ldots, L,\notag\\
  \mbox{Flow model:}\quad\frac{\partial u}{\partial t} &= \mathcal{F}(u,t),\label{LaDA_system}
\end{align}
where $L$ is the number of tracers, and $\dot{W}_l$ is a white noise, representing the observational noise. Each $x_l$ is a two-dimensional vector containing the displacements of the $l$-th tracer, and $u$ is the two-dimensional velocity field (e.g., the surface of the ocean). Notably, the flow velocity in the observational process is typically a highly nonlinear function of displacement, making the LaDA a challenging nonlinear problem.

\subsection*{UQ in LaDA}
One important UQ topic in LaDA is quantifying the uncertainty reduction in the estimated flow field as a function of the number of tracers $L$. The comparison is between the posterior distribution from LaDA and the prior one. The latter is the statistical equilibrium solution (i.e., the attractor) of the chaotic flow model. It is the best inference of the flow field in the absence of observations. With two distributions involved in the comparison, the relative entropy \eqref{Relative_Entropy} is a natural choice to measure the uncertainty reduction.

\begin{figure}[htbp]
	\begin{center}
		\hspace*{-0cm}\includegraphics[width=8cm]{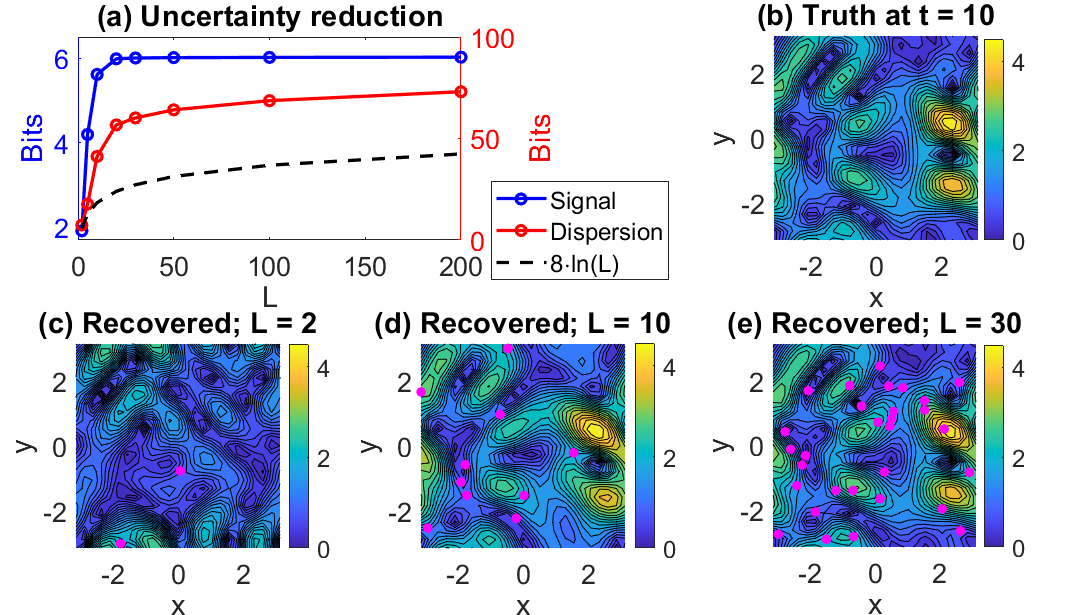}
	\end{center}
	\caption{UQ in LaDA. Panel (a): uncertainty reduction in the signal and dispersion parts as a function of $L$. Panels (b)--(e): comparing the true flow field with the recovered ones using 2, 10 and 50 tracers.   }
	\label{tracers_information}
\end{figure}

Under certain conditions \cite{chen2014information}, closed analytic formulae are available for finding the posterior distribution of the nonlinear LaDA problem, which allows for using rigorous analysis to quantify the uncertainty reduction for the long-term behavior. By applying \eqref{Relative_Entropy_Gaussian}, it can be shown that the signal part of the relative entropy converges to a constant as $L\to\infty$, where the posterior mean $\mu$ in \eqref{Relative_Entropy_Gaussian} is replaced by the true flow field. This is intuitive as the tracers become dense in the domain, and the observational error cancels off at each location when $L\to\infty$. On the other hand, the dispersion part in \eqref{Relative_Entropy_Gaussian} will never converge. Instead, it grows at $\ln (L)$. Although it confirms that deploying extra tracers will always bring additional information, the increment of the uncertainty reduction will decrease as $L$ becomes large. Specifically, the logarithm dependence on $L$ implies that reducing the uncertainty by a fixed amount requires an exponential increase in the number of tracers. This means, in practice, once the information gain reaches a certain level, there is no need to continue deploying tracers as the uncertainty reduction from the additional tracers becomes marginal. Figure \ref{tracers_information} shows a numerical simulation that validates the theoretic conclusion.

\section*{Role of the Uncertainty in Diagnostics}
\subsection*{Parameter estimation}
Consider the following two-dimensional model
\begin{equation}\label{Linear_Oscillator}
  \frac{\d x}{\d t} = ay,\qquad \frac{\d y}{\d t} = bx,
\end{equation}
which is a linear system with respect to both the parameter $\theta=(a,b)^\mathtt{T}$ and the state variable $(x,y)^\mathtt{T}$. Assume the observational data of $x$ and $y$ and their time derivatives $\dot{x}:=\d x/\d t$ and $\dot{y}:=\d y/\d t$ are available at time $t_i$ for $i=1,\ldots, I$. Define a matrix $M_i=
                                                                                                                                                  \begin{pmatrix}
                                                                                                                                                    y_i & 0 \\
                                                                                                                                                    0 & x_i \\
                                                                                                                                                  \end{pmatrix}
                                                                                                                                                $,
which takes values at time $t_i$. To estimate the parameters $\theta$, linear regression can be easily applied to find the least-squares solution,
\begin{equation}\label{regression_parameter_estimation}
  \theta = \left(\sum_i M_i^\mathtt{T}M_i\right)^{-1}\left(\sum_i M_i^\mathtt{T}z_i\right),
\end{equation}
where $z_i = \left(\dot{x}_i, \dot{y}_i\right)^\mathtt{T}$. Writing down the component-wise form of \eqref{regression_parameter_estimation} yields
\begin{equation}\label{regression_parameter_estimation_a}
  a = \left(\sum_i y_i^2\right)^{-1}\left(\sum_i y_i\dot{x}_i\right).
\end{equation}

Now consider a different situation, where only $x_i$ and $\dot{x}_i$ are observed. This is typical in many applications where only partial observations are available. Since estimating the parameter $a$ requires the information of $y$ as shown in \eqref{regression_parameter_estimation_a}, the state of $y$ at each time $t_i$ needs to be estimated, which naturally introduces uncertainty. Assume the estimated state of $y_i$ is given by a Gaussian distribution, and $y_i$ is written into the Reynolds decomposition form as $y_i = \langle y_i\rangle + y_i^\prime$. The parameter estimation problem can be regarded as repeatedly drawing samples from the distribution of $y_i$, plugging them into \eqref{regression_parameter_estimation_a}, and then taking the average for the evaluation of the terms involving $y_i$ to reach the maximum likelihood solution. In light of the fact from \eqref{mean_fluctuation_variable} that $\langle y_i^2\rangle=\langle y_i\rangle^2 + \langle(y_i^\prime)^2\rangle$, the above procedure essentially gives the following modified version of \eqref{regression_parameter_estimation_a},
\begin{equation}\label{regression_parameter_estimation_a_UQ}
  a = \left(\sum_i \left(\langle y_i\rangle^2 + \langle(y_i^\prime)^2\rangle\right)\right)^{-1}\left(\sum_i \langle y_i\rangle\dot{x}_i\right),
\end{equation}
where the term $\langle(y_i^\prime)^2\rangle$ arrives due to the average of the nonlinear function $y_i^2$ in \eqref{regression_parameter_estimation_a}. This is somewhat surprising at a glance since the underlying model \eqref{Linear_Oscillator} is linear. Nevertheless, the formula in \eqref{regression_parameter_estimation_a_UQ} is intuitive. Consider estimating $a$ from $\dot{x} = a y_i$ with data at only one time instant $t_i$. If $y_i\neq0$ is deterministic, then $a=\dot{x}/y_i$. If $y_i$ contains uncertainty and satisfies a Gaussian distribution, then the reciprocal distribution $\dot{x}/y_i$ is no longer Gaussian. The variance of $y_i$ affects the mean value $\langle\dot{x}/y_i\rangle$, which is different from $\dot{x}/\langle y_i\rangle$. Following the same logic, as $y_i^2$ on the right-hand side of the least-squares solution \eqref{regression_parameter_estimation_a} is nonlinear, the additional term $\langle(y_i^\prime)^2\rangle$ appears in \eqref{regression_parameter_estimation_a}, which affects the parameter estimation skill. Therefore, uncertainty can play a significant role even in diagnosing a linear system where nonlinearity appears in the diagnostics (parameter estimation formula).

As a numerical illustration, assume enough data is generated from \eqref{Linear_Oscillator} with the true parameters $a=2$ and $b=-2$. When the uncertainty, namely the variance of estimating $y_i$, is given by $\langle(y_i^\prime)^2\rangle=0.5,1$ and $2$ at all times $t_i$, the estimated parameter of $a$ becomes $a=1.000$, $0.672$, and $0.404$, respectively. As the estimated parameter becomes less accurate, the residual term in the regression increases, accounting for the effect of the input uncertainty.

In practice, expectation-maximization (EM) iterative algorithms can be applied to alternatively update the estimated parameter values and recover the unobserved states of $y$ \cite{chen2023stochastic}. In the E step, the state of $y$ can be inferred via DA, using the current parameter values $\theta$. In the M step, the maximum likelihood estimate is applied to infer the parameters given the estimated state of $y$. Likewise, data augmentation can be used to sample the trajectories of the unobserved variables. This can be incorporated into other parameter estimation methods, such as the Markov Chain Monte Carlo \cite{papaspiliopoulos2013data}.

\subsection*{Eddy identification}
Let us now explore how simple UQ tools advance the study of realistic problems. Oceanic eddies are dynamic rotating structures in the ocean. The primary goal is to show that nonlinearity in the eddy diagnostic will make UQ play a crucial role. Providing a rigorous definition of an eddy and discussing the pros and cons of different eddy identification methods are not the main focus here. Mesoscale eddies are major drivers of the transport of momentum, heat, and mass, as well as biochemical and biomass transport and production in the ocean. The study of ocean eddies is increasingly important due to climate change and the vital role eddies play in the rapidly changing polar regions and the global climate system.

Due to the complex spatial and dynamical structure of eddies there is no universal criterion for eddy identification. The Okubo-Weiss (OW) parameter \cites{okubo1970horizontal, weiss1991dynamics} is a widely used approach based on physical properties of the ocean flow:
    \begin{equation}\label{OW_Parameter}
        \operatorname{OW} = s_\mathrm{n}^2 + s_\mathrm{s}^2 - \omega^2,
    \end{equation}
    where the normal strain, the shear strain, and the relative vorticity are given by
    \begin{equation*}
        s_\mathrm{n} = u_x-v_y, \quad s_\mathrm{s} = v_x+u_y, \quad\mbox{and}\quad\omega = v_x-u_y,
    \end{equation*}
    respectively, with $\mathbf{u}=(u,v)$ the two-dimensional velocity field and the shorthand notation $u_x:=\partial u/\partial x$. When the OW parameter is negative, the relative vorticity is larger than the strain components, indicating vortical flow. The OW parameter is an Eulerian quantity based solely on a snapshot of the ocean velocity field. There are other quantities to identify eddies, such as the Lagrangian descriptor \cite{vortmeyer2016detecting}, based on a sequence of snapshots.

The use of these eddy identification diagnostics requires knowing the exact flow field. However, uncertainties may appear in state estimation, for example, in the marginal ice zone when the ocean is estimated via the LaDA using a limited number of ice floe trajectories. Eddy identification exploiting the mean estimate of the flow field can lead to large biases.

The contribution of the uncertainty in affecting the OW parameter can be seen by applying the Reynolds decomposition to each component in \eqref{OW_Parameter}, e.g., $u_x = \bar{u}_x + u_x'$. By sampling multiple realizations from the posterior distribution of the estimated flow field, the expectation of the OW parameter values applying to each sampled flow field is given by
\begin{equation}\label{OW_mean}
\begin{split}
\mathbb{E}[\mathrm{OW}(\mathbf{u})] =\mathrm{OW}(\bar{\mathbf{u}})&+ \langle(u_x')^2\rangle - 2\langle u_x'v_y'\rangle\\
        & + \langle(v_y')^2\rangle + 4\langle v_x'u_y'\rangle.
\end{split}
\end{equation}
On the right-hand side of \eqref{OW_mean}, $\mathrm{OW}(\bar{\mathbf{u}})$ is the OW parameter applying to the estimated mean flow field. The remaining terms are all nonlinear functions of the fluctuation part due to the uncertainty in state estimation. In general, as most of the eddy identification criteria are nonlinear with respect to the flow field as in \eqref{OW_Parameter}, changing the order of taking the expectation and applying the eddy diagnostic criterion will significantly impact the results. Notably, even though \eqref{regression_parameter_estimation_a_UQ} and \eqref{OW_mean} are in entirely different contexts, they share the same essence that uncertainty plays a significant role when nonlinear terms appear in the diagnostics.

\begin{figure}[ht]
	\begin{center}
		\hspace*{-0cm}\includegraphics[width=8cm]{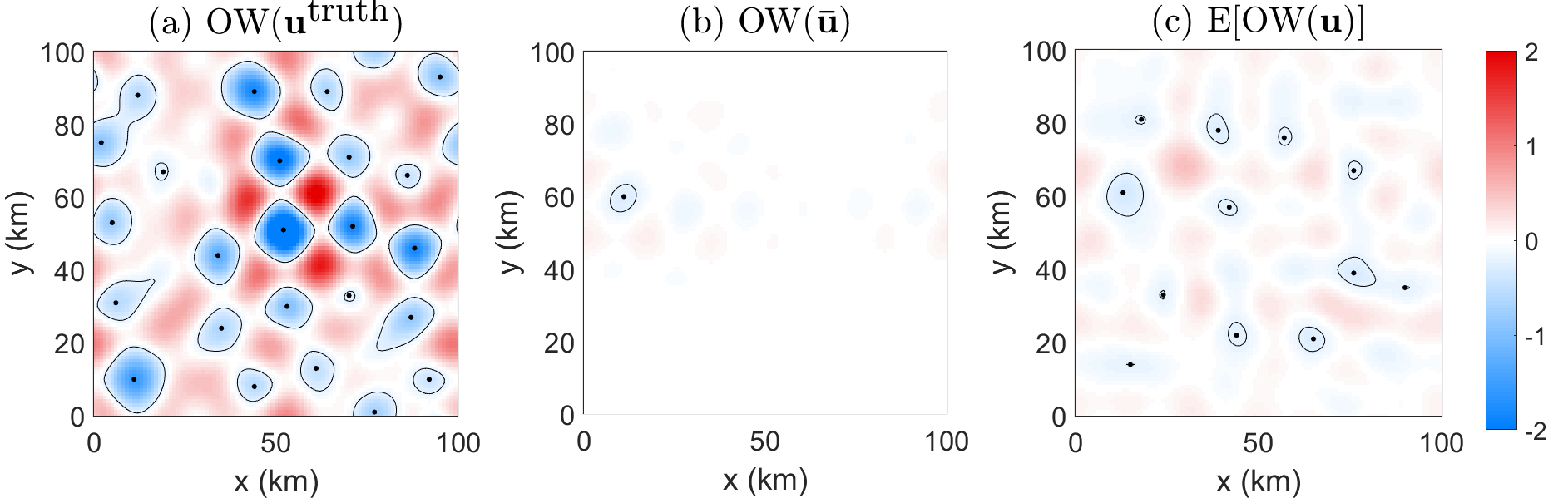}
	\end{center}
	\caption{Eddy identification using the OW parameter. Panel (a): the OW parameter based on the true flow field. Panel (b): the OW parameter based on the estimated mean flow field from LaDA using a small number of the observed tracers. Panel (c): the expected value of the OW parameter based on multiple flow fields sampled from the posterior distribution.   }
	\label{owmeancomparison2}
\end{figure}

Panels (b)--(c) of Figure \ref{owmeancomparison2} numerically illustrate the difference between $\mathbb{E}[\mathrm{OW}(\mathbf{u})]$ and $\mathrm{OW}(\bar{\mathbf{u}})$, where the flow field is inferred from LaDA using a small number of the observed tracers. In the presence of such a large uncertainty, the main contribution of the estimated flow field comes from the fluctuation part, while the mean estimation has weak amplitudes. Consequently, very few eddies are identified using $\mathrm{OW}(\bar{\mathbf{u}})$. The expected value $\mathbb{E}[\mathrm{OW}(\mathbf{u})]$ improves the diagnostic results. However, it is still far from the truth in Panel (a) due to the large uncertainty of $\mathrm{OW}(\mathbf{u})$ computed by using different realizations of $\mathbf{u}$ sampled from the posterior distribution of the LaDA. This is again similar to the intrinsic inaccuracy in the parameter estimation problem \eqref{regression_parameter_estimation_a_UQ} when uncertainty appears. Notably, collecting $\mathrm{OW}(\mathbf{u})$ associated with different sampled flow realizations allows us to compute the PDF $p(\mathrm{OW}(\mathbf{u}))$, which provides a probabilistic view of eddy identification. Such a probabilistic eddy identification framework allows for assigning a probability of the occurrence of each eddy and the PDFs of lifetime and size of each eddy \cite{covington2024probabilistic}.

\section*{UQ in Advancing Efficient Modeling}
The governing equations of many complex turbulent systems are given by nonlinear partial differential equations (PDEs). The spectral method remains one of the primary schemes for finding numerical solutions. For simplicity, consider a scalar field $u(x,t)$, which satisfies a PDE with quadratic nonlinearities, as in many fluids and geophysical problems. Assume Fourier basis functions $\exp(ikx)$ are utilized, where $k$ is the wavenumber, and $x$ is the spatial coordinate. Denote by $\hat{u}_k(t)$ the time series of the spectral mode associated with wavenumber $k$. Apply a finite-dimensional truncation to retain only the modes $k\in[-K,K]$. The resulting equation of $\hat{u}_k(t)$ typically has the following form,
\begin{equation}\label{Deterministic_model_one_mode}
\begin{split}
    \frac{\d\hat{u}_k(t)}{\d t} =& \big(-d_k +i\omega_k\big)\hat{u}_k(t)+{\hat{f}_k(t)}\\&\quad+\sum_{-K\leq m\leq K} c_{m,k}\hat{u}_m\hat{u}_{k-m},
\end{split}
\end{equation}
where $K\gg 1$ represents the resolution of the model, $i$ is the imaginary unit, $d_k$ and $\omega_k$ are both real, while $\hat{f}_k(t)$ is complex.
The first term on the right-hand side of \eqref{Deterministic_model_one_mode} is linear, representing the effect of damping/dissipation $-d_k<0$ and phase $\omega_k$. The second term is deterministic forcing.
The last term sums over all the quadratic nonlinear interactions projected to mode $k$. Typically, a large number of such nonlinear terms will appear in the summation, which is one of the main computational costs in solving \eqref{Deterministic_model_one_mode}.

Due to the turbulence of these systems, many practical tasks, such as the statistical forecast and DA, require obtaining a forecast distribution by repeatedly running the governing equation. Each run is already quite costly, so the ensemble forecast is usually computationally prohibitive. Therefore, developing an appropriate stochastic surrogate model is desirable, aiming to significantly accelerate computational efficiency. Stochasticity can mimic many of the features in turbulent systems. Since the goal is to reach the forecast statistics instead of a single individual forecast trajectory, once appropriate UQ is applied to guide the development of the stochastic surrogate model, it can reproduce the statistical features of the underlying nonlinear deterministic system.

In many applications, a large portion of the energy is explained by only a small number of the Fourier modes. Yet, the entire set of the Fourier modes has to be solved together in the direct numerical simulation to guarantee numerical stability. Therefore, reducing the computational cost by using stochastic surrogate models is two-folded. First, in the governing equation of each mode, the heavy computational burden of calculating the summation of a large number of nonlinear terms will be replaced by computing a few much cheaper stochastic terms. Second, as the governing equation of each mode becomes independent when the nonlinear coupling is replaced by the stochastic terms, only the leading few Fourier modes need to be retained in such an approximate stochastic system, saving a large amount of computational storage. This also allows a larger numerical integration time step since stiffness usually comes from the governing equations of the small-scale modes. One simple stochastic model is to replace all the nonlinear terms in the governing equation of each mode by a single stochastic term $\dot{W}$ representing white noise, which works well if the statistics of the mode are nearly Gaussian. The resulting stochastic model reads \cite{gardiner2004handbook}:
 \begin{equation}\label{Stochastic_model_one_mode}
    \frac{\d\hat{u}_k(t)}{\d t} = \big(-d_k +i\omega_k\big)\hat{u}_k(t)+{\hat{f}_k(t)}+\sigma_k\dot{W}_k.
\end{equation}
For systems with strong non-Gaussian statistics, other systematic methods can be applied to develop stochastic surrogate models \cite{chen2023stochastic}.

UQ plays a crucial role in calibrating the stochastic surrogate model. Specifically, the parameters in the stochastic surrogate model are optimized so that the two models have the same forecast uncertainty, which is crucial for DA and ensemble prediction. This is often achieved by matching a set of key statistics of the two models, especially the equilibrium PDFs and the decorrelation time. For \eqref{Stochastic_model_one_mode}, simple closed analytic formulae are available for model calibration in reproducing the forecast uncertainty, allowing it to be widely used in many practical problems. See the online supplementary document for more details.

\begin{figure}[ht]
	\begin{center}
		\hspace*{-0cm}\includegraphics[width=8cm]{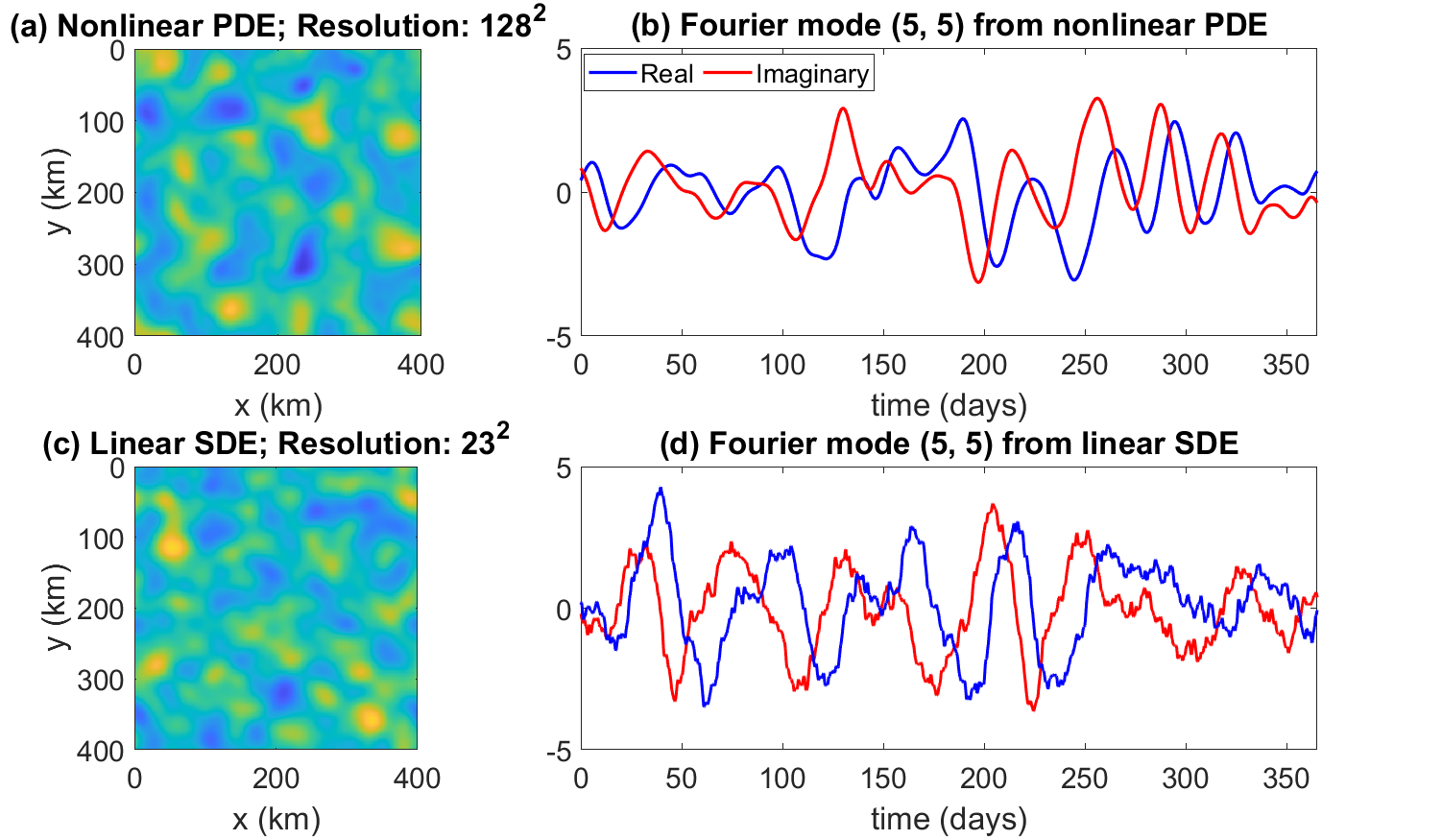}
	\end{center}
	\caption{Comparison of the simulations from a nonlinear PDE model and a stochastic model, showing the latter can reproduce the forecast statistics of the former. Panels (a) and (c): snapshots of stream functions of the upper-layer ocean. Panels (b) and (d): time series of mode $(5,5)$. These are two different realizations, so there is no point-wise correspondence between them in the snapshots or time series.  }
	\label{SDE_PDE_Approx}
\end{figure}

Panels (a) of Figure \ref{SDE_PDE_Approx} shows a snapshot of the stream function from a two-layer quasi-geostrophic (QG) model, which is given by a set of nonlinear PDEs with the spatial resolution of each layer being $128\times128$. In comparison, Panel (c) shows a snapshot of the spatial field by running a set of stochastic surrogate models, each calibrated by capturing the forecast uncertainty of one Fourier mode of the two-layer QG model \cite{covington2022bridging}. Although the spatial resolution is only $23\times 23$, which contains less than 2\% of the total number of modes, the spatial field is qualitatively similar to the QG model. Panels (b) and (d) show the time series of Fourier mode $(5,5)$ associated with the two models. They behave similarly as they have the same temporal correlation and PDFs. An additional example of using the stochastic model \eqref{Stochastic_model_one_mode} to approximate the KdV equation can be found in \cite{majda2019statistical}.

\section*{Conclusions}
This paper exploits simple examples to provide the basic concepts of UQ. It does not aim to provide a single, comprehensive definition, as various communities currently have different interpretations of UQ. Nevertheless, the information quantities \eqref{Shannon_Entropy} and \eqref{Relative_Entropy} are natural measurements to quantify the uncertainty. As a developing field within applied mathematics and interdisciplinary research, UQ continues to evolve. One of the goals of this paper is to offer new insights into how these ideas can be applied across different fields, helping to reveal the commonalities and practical advantages of diverse approaches. As uncertainty is ubiquitous, incorporating UQ into analysis or strategic planning is essential to facilitate understanding nature and problem-solving in almost all disciplines. UQ has become an important tool in tackling observational data \cite{braverman2021post}, forecasting turbulent systems \cite{majda2012lessons}, and in studying inverse problems \cite{dashti2011uncertainty}. UQ is not a stand-alone subject. It strongly relates to and advances the development of new techniques for mathematical reasoning, computational algorithms, effective modeling, data analysis, machine learning, and various applications. The areas in science, engineering, and technology requiring careful UQ are rapidly expanding and demands new ideas for coping with challenging issues.
\medskip

\noindent
{\bf Acknowledgements.} N.C. and M.A. are grateful to acknowledge the support of the Office of Naval Research (ONR) N00014-24-1-2244 and Army Research Office (ARO) W911NF-23-1-0118. S.W. acknowledges the financial support provided by the EPSRC Grant No. EP/P021123/1 and the support of the William R. Davis '68 Chair in the Department of Mathematics at the United States Naval Academy. The authors thank Dr. Jeffrey Covington for helping with some of the figures.\medskip

\noindent
{\bf Supplementary Document and Codes Availability.} The arXiv version of this article includes a supplementary document (\url{https://arxiv.org/abs/2408.01823}), which contains a comprehensive tutorial of many numerical examples shown in this article and beyond. The codes for these examples, written in both MATLAB and Python, are available from GitHub at \url{https://github.com/marandmath/UQ_tutorial_code}.
\bibliography{references}

\newpage

\onecolumn

\noindent\textbf{\Large Supplementary Document}

\tableofcontents

\vspace{0.5cm}

\noindent The notations in this supplementary document are slightly different from those in the main text. Throughout this supplementary document, \textbf{boldface} letters will be will exclusively used to denote vectors for the sake of mathematical clarity. In this regard, we will specifically use \textbf{l}owercase boldface letters to denote column vectors, while we will use \textbf{U}ppercase boldface letters to denote matrices.

\section{Basic Probabilities}
Consider a scalar random variable $X$ with PDF $p(x)$. Denote by $E[\cdot]$ the expectation of a random variable. The $n$-th moment and central moment are defined respectively as,
\begin{equation}\label{Moments_SI}
\begin{split}
	\mu_n &= E[X^n] = \int_{-\infty}^{\infty} x^np(x)\d x,\qquad \mbox{for~} n\geq1,\\
	\tilde{\mu}_n &= E[(X-\mu)^n] = \int_{-\infty}^{\infty} (x-\mu)^np(x)\d x,\qquad\mbox{for~} n\geq2.
\end{split}
\end{equation}
The expressions of the moments are sometimes denoted as $\mu^n = \langle X^n \rangle$ for notational simplicity, where $\langle\boldsymbol{\cdot}\rangle:=\int_{-\infty}^\infty\boldsymbol{\cdot}\, p(x)\d x$ stands for the statistical average over the state variable. The two most widely used moments are the mean and the variance.
\begin{itemize}
  \item The mean $\mu$ is the first moment, representing the average value of the random variable,
  \item The variance $\sigma^2:=\tilde\mu_2$ is the second central moment, measuring how far a set of numbers spreads out from the mean.
\end{itemize}
The square root of the variance $\sigma$ is named the standard deviation, which has the same unit as the random variable. For multi-dimensional random variables, the second central moment is called the covariance matrix, which measures the joint variability between the components of the random variable. For $n\geq 3$, the standardized moments provide more intricate characterizations of a random variable and the associated PDF. A standardized moment is a normalized central moment (with $n\geq3$). Normalization is typically a division by an expression of the standard deviation, which renders the moment scale invariant. The third $(n=3)$ and the fourth ($n=4$) standardized moments are named skewness and kurtosis, which characterize the asymmetry about the mean and the ``tailedness'' associated with the PDF of a real-valued random variable, respectively. They are expressed as,
\begin{equation}\label{Skewness_Kurtosis_SI}
  \mbox{Skewness }= E\left[\left(\frac{X-\mu}{\sigma}\right)^3\right] = \frac{\tilde\mu_3}{\tilde\mu_2^{3/2}}\qquad \mbox{and}\qquad \mbox{Kurtosis } = E\left[\left(\frac{X-\mu}{\sigma}\right)^4\right] = \frac{\tilde\mu_4}{\tilde\mu_2^{2}}.
\end{equation}
When the PDF is symmetric, the skewness is zero. A positive (negative) skewness commonly indicates that the tail is on the right (left) side of the distribution. A PDF with kurtosis being the standard value of $3$ has the same tail behavior as the Gaussian distribution. When the kurtosis is smaller (larger) than $3$, the random variable produces fewer (more) extreme outliers than the corresponding Gaussian counterpart. Therefore, the associated PDF has lighter (heavier) tails than the Gaussian distribution. Excess kurtosis is sometimes used, defined as kurtosis minus $3$. Skewness and kurtosis are important indicators for measuring the non-Gaussianity of a PDF.

\section{Computing Shannon's Entropy and Relative Entropy}
\subsection{Shannon's entropy}
{\color{blue}{
The codes that correspond to the content of this subsection are:
\begin{itemize}
  \item MATLAB code: \textbf{Computing\_Entropy.m}
  \item Python code: \textbf{Computing\_Entropy.py}
\end{itemize}}}

\noindent Denote by $p(x)$ the PDF of a random variable $X$. Recall the definition of Shannon's entropy:
\begin{equation}\label{Shannon_Entropy_SI}
    \mathcal{S}(p) = -\int p(x) \ln(p(x)) \d x.
\end{equation}
For certain distributions, $\mathcal{S}(p)$ can be written down explicitly.

\paragraph{Gaussian distribution.} If $p\sim\mathcal{N}_m(\boldsymbol{\mu}, \mathbf{R})$ is an $m$-dimension Gaussian distribution, where $\boldsymbol{\mu}$ and $\mathbf{R}$ are the mean vector and the covariance matrix, respectively, then the Shannon entropy has the following explicit form:
  \begin{equation}\label{Shannon_Entropy_Gaussian_SI}
    \mathcal{S}(p) = \frac{m}{2}(1+\ln2\pi) + \frac{1}{2}\ln\det(\mathbf{R}).
  \end{equation}
The explicit formula in \eqref{Shannon_Entropy_Gaussian_SI}, when $m=1$, can be derived as follows. The one-dimensional Gaussian is given by,
\begin{equation}\label{One_D_Gaussian_SI}
    p(x) = \frac{1}{\sqrt{2\pi R}}\exp\left(-\frac{(x-\mu)^2}{2R}\right).
\end{equation}
Plug \eqref{One_D_Gaussian_SI} into \eqref{Shannon_Entropy_SI} and replace only the $p(x)$ in the logarithm function but keep the abstract form of the other $p(x)$ that represents the weight.
Then the Shannon entropy becomes,
\begin{equation*}
\begin{split}
    -\int p(x)\ln p(x)\d x &= -\int p(x) \left( -\frac{1}{2}\ln \left( 2\pi R\right) - \frac{(x-\mu)^2}{2R}\right)\d x\\
    & = \frac{1}{2}\ln (2\pi R) + \frac{1}{2R} \int(x-\mu)^2 p(x)\d x\\
    & = \frac{1}{2}\ln (2\pi) + \frac{1}{2}\ln R + \frac{1}{2},
\end{split}
\end{equation*}
where the first term in the second equality comes from the fact that $\int p(x)\d x=1$, while the third term in the third equality uses the definition of the variance $R = \int(x-\mu)^2 p(x)\d x$.

\paragraph{Gamma distribution.} Gamma distribution is a one-dimensional two-parameter family of continuous probability distributions. The exponential distribution, Erlang distribution, and chi-squared distribution are special cases of the gamma distribution. There are two equivalent parameterizations for Gamma distribution: (a) with a shape parameter $k$ and a scale parameter $\theta$, and (b) with a shape parameter $\alpha=k$ and inverse scale parameter $\beta = 1/\theta$, known as the rate parameter. In each of these forms, both parameters are positive real numbers. We will focus on the first form in this document.

The PDF of the gamma distribution is given by,
\begin{equation}\label{Gamma_distribution_SI}
  p(x) = \frac{1}{\Gamma(k)\theta^k}x^{k-1}e^{-\frac{x}{\theta}},
\end{equation}
where,
\begin{equation}\label{Gamma_function_SI}
  \Gamma(k) = \int_0^\infty t^{k-1}e^{-t}\d t,\qquad\mbox{for } k>0.
\end{equation}
Different from the Gaussian distribution that is uniquely determined by the mean and the covariance (the first two centralized moments), the Gamma distribution has non-vanishing higher order moments. The leading four moments can be expressed in terms of the shape and scale parameters as,
\begin{equation}\label{Gamma_distribution_stats_SI}
  \mbox{mean}=k\theta, \qquad \mbox{variance}=k\theta^2,\qquad\mbox{skewness}=\frac{2}{\sqrt{k}},\qquad\mbox{and}\qquad\mbox{kurtosis}=\frac{6}{k}.
\end{equation}
Shannon's entropy also has an explicit expression for the case of a Gamma distribution, which is,
\begin{equation}\label{Gamma_distribution_entropy_SI}
  \mathcal{S} = k + \ln \theta + \ln\Gamma(k) + (1-k)\psi(k),
\end{equation}
where $\psi(k)$ is the digamma function defined by,
\begin{equation}\label{Digamma_function_SI}
  \psi(k) = \frac{\d }{\d x}\ln\Gamma(k) = \frac{\Gamma^\prime(k)}{\Gamma(k)}.
\end{equation}

\paragraph{General PDFs.} There is often no closed-form formula for computing Shannon's entropy \eqref{Shannon_Entropy_SI} for a general distribution. Nevertheless, numerical methods, such as the trapezoidal rule, can be applied, at least when the dimension of the state variable $x$ remains relatively low. Note that $p(x)$ usually covers an infinite domain while the numerical approximation of $p(x)$ is generally based on a finite interval. Therefore, normalization is often required for $p(x)$ to guarantee the integration of $p(x)$ within the finite interval for the numerical solution is one before plugging it into \eqref{Shannon_Entropy_SI} to compute Shannon's entropy.

\subsection{Relative entropy}
{\color{blue}{
The codes that correspond to the content of this subsection are:
\begin{itemize}
  \item MATLAB code: \textbf{Computing\_Relative\_Entropy.m}
  \item Python code: \textbf{Computing\_Relative\_Entropy.py}
\end{itemize}}}

\noindent Recall the definition of the relative entropy:
\begin{equation}\label{Relative_Entropy_SI}
    \mathcal{P}(p,p^M)= \int p(x) \ln\left(\frac{p(x)}{p^M(x)}\right) \d x.
  \end{equation}
When both $p\sim\mathcal{N}_m(\boldsymbol{\mu},\mathbf{R})$ and $p^M\sim\mathcal{N}_m(\boldsymbol{\mu}^M,\mathbf{R}^M)$ are $m$-dimensional Gaussians, the relative entropy has the following explicit formula known as the signal-dispersion decomposition,
\begin{equation}\label{Relative_Entropy_Gaussian_SI}
    \mathcal{P}(p,q)= \underbrace{\frac{1}{2}\Big[(\boldsymbol{\mu}-\boldsymbol{\mu}^M)^\mathtt{T}(\mathbf{R}^M)^{-1}(\boldsymbol{\mu}-\boldsymbol{\mu}^M)\Big]}_{\mbox{Signal}} + \underbrace{\frac{1}{2}\Big[\mbox{tr}(\mathbf{R}(\mathbf{R}^M)^{-1})-m-\ln\mbox{det}(\mathbf{R}(\mathbf{R}^M)^{-1})\Big]}_{\mbox{Dispersion}},
\end{equation}
where 'tr' and 'det' are the trace and determinant of a matrix, respectively. The first term on the right-hand side of \eqref{Relative_Entropy_Gaussian_SI} is called `signal', which measures the lack of information in the mean weighted by the model covariance. The second term involving the covariance ratio is called `dispersion'.

For two general distributions, computing the relative entropy has to go through the formula in \eqref{Relative_Entropy_SI}. In practice, the PDFs are usually estimated from data using histogram or kernel density estimations. Since the number of data points used to estimate PDFs is finite, the tails of PDFs are typically underestimated. Sometimes, the $p^M(x)$ is even zero for certain values of $x$ close to the tails within the interval range. However, $p^M(x)$ appears in the denominator. These small numerical errors in estimating the PDFs may significantly change the value of relative entropy. Therefore, numerical remedies are needed. One simple remedy is applying the following two-step procedure to the estimated PDF of $p^M(x)$ (as well as $p(x)$):
\begin{itemize}
  \item Step 1. Set up a threshold $\epsilon$, for example $\epsilon=10^{-5}$. Set $p^M(x) = \epsilon$ for those $x$ such that the originally estimated $p^M(x)< \epsilon$. The resulting PDF is denoted by $\tilde{p}^M(x)$.
  \item Step 2. Normalize $\tilde{p}^M(x)$ by dividing it by $\int p^M(x)\d x$. Set the result to be $p^M(x)$.
\end{itemize}

\begin{figure}[ht!]
	\begin{center}
		\hspace*{-0.5cm}\includegraphics[width=17.5cm]{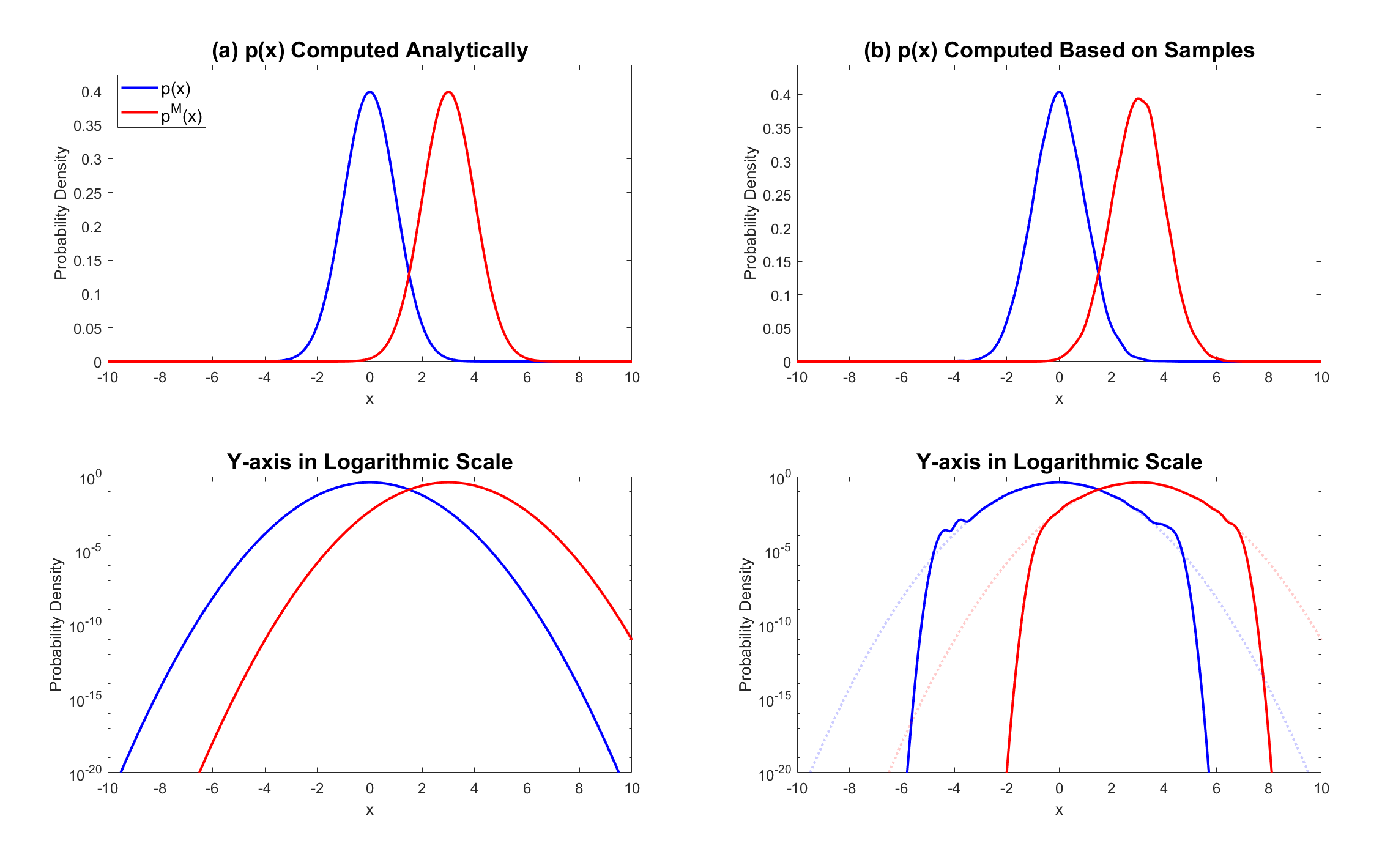}
	\end{center}
	\caption{Comparison of the Gaussian distributions constructed in different ways. Column (a): Given the mean and variance, each PDF is built using the analytic formula \eqref{One_D_Gaussian_SI}. Column (b): Given the mean and variance, $10000$ samples are drawn from the corresponding Gaussian distributions. The PDFs shown here are then built upon these samples. The second row shows the same PDFs as the first, but the y-axis is plotted using a logarithmic scale to better illustrate the tail behavior. Note that the Gaussian PDFs have parabolic profiles when shown in a logarithmic scale due to the quadratic term in the exponential. For reference, the dotted curves in the last panel represent the parabolas.   }
	\label{SI_Relative_Entropy_Numerics}
\end{figure}

Figure \ref{SI_Relative_Entropy_Numerics} compares the Gaussian distributions constructed in different ways. In Column (a), each PDF is built using the analytic formula given the mean and variance \eqref{One_D_Gaussian_SI}. Due to the analytic expressions of these PDFs, the calculated relative entropy using the direct definition by taking the numerical integration \eqref{Relative_Entropy_SI} has the same result as that using the explicit formula \eqref{Relative_Entropy_Gaussian_SI}. In Column (b), we consider a situation in which the exact expression of the PDFs is unknown. The available information is a finite number of sample points. We adopt a relatively large number of samples --- $10000$. Assume these samples are taken from a Gaussian distribution, which is unknown in practice. Then, a kernel density estimation is applied to characterize the PDF numerically. It is seen in Column (b) that the core of the PDF is approximated quite well. However, the tails of the distribution, corresponding to the extreme events in the samples, are severely underestimated. Therefore, if the direct definition \eqref{Relative_Entropy_SI} is utilized to compute the relative entropy, the result will be infinite since the tail probability of $p^M$ becomes zero when it is far from the center but within the domain $[-10,10]$. Therefore, the above remedy becomes essential, leading to a value similar to the truth.

\section{UQ in Dynamical Systems}
\subsection{Uncertainty propagation in the linear damped system}
{\color{blue}{
The codes that correspond to the content of this subsection are:
\begin{itemize}
  \item MATLAB code: \textbf{Linear\_System\_UQ.m}
  \item Python code: \textbf{Linear\_System\_UQ.py}
\end{itemize}}}
\noindent The code will reproduce Figure \ref{ODE_UQ} in the main text.
\subsection{Uncertainty propagation in the chaotic Lorenz 63 model}
{\color{blue}{
The codes that correspond to the content of this subsection are:
\begin{itemize}
  \item MATLAB code: \textbf{L63\_UQ.m}
  \item Python code: \textbf{L63\_UQ.py}
\end{itemize}}}
\noindent The code will reproduce Figure \ref{L63_UQ} in the main text.

\section{Uncertainty Reduction via Bayes Formula and Data Assimilation (DA)}
\subsection{Uncertainties in posterior distributions}
{\color{blue}{
The codes that correspond to the content of this subsection are:
\begin{itemize}
  \item MATLAB code: \textbf{Bayes\_Formula.m}
  \item Python code: \textbf{Bayes\_Formula.py}
\end{itemize}}}

\noindent Recall Bayes' formula:
\begin{equation}\label{Bayes_formula_SI}
{ \underbrace{p(\mathbf{u}|\mathbf{v})}_{\mbox{posterior}}} \propto~ { \underbrace{p(\mathbf{u})}_{\mbox{prior}}}~\times\underbrace{p(\mathbf{v}|\mathbf{u})}_{\mbox{likelihood}}.
\end{equation}
The observation or measurement is given by,
\begin{equation}\label{Observations_Bayes_SI}
{\mathbf{v}}=\mathbf{G}\mathbf{u}+\boldsymbol{\epsilon},
\end{equation}
where $\mathbf{G}$ is the observational operator and $\boldsymbol{\epsilon}$ is a zero-mean Gaussian white noise. Denote with $p(\mathbf{u})\sim\mathcal{N}(\boldsymbol{\mu}_f,\mathbf{R}_f)$ the prior distribution. Further denote by $\mathbf{R}^o$ the observational noise, namely the covariance matrix of each measurement, $\mathbf{v}$. Since the prior and the likelihood are both Gaussian, the posterior distribution $p(\mathbf{u}|{\mathbf{v}})\sim\mathcal{N}(\boldsymbol{\mu}_a,\mathbf{R}_a)$ in \eqref{Bayes_formula_SI} can be written in a more explicit way as,
\begin{equation}\label{Posterior_SI}
  p(\mathbf{u}|\mathbf{v})\propto p(\mathbf{u})p(\mathbf{v}|\mathbf{u}) = e^{-\frac{1}{2}J(\mathbf{u})},
\end{equation}
where the quadratic term in the exponential is given by,
\begin{equation}\label{Posterior_Supp_SI}
  J(u) = (\mathbf{u}-\boldsymbol{\mu}_{f})^*\left(\mathbf{R}_f\right)^{-1}(\mathbf{u}-\boldsymbol{\mu}_{f}) + (\mathbf{v}-\mathbf{g}\mathbf{u})^*\left(\mathbf{R}^o\right)^{-1}(\mathbf{v}-\mathbf{g}\mathbf{u}).
\end{equation}
with $\cdot^*$ being the transpose (or conjugate transpose if $\mathbf{u}$ is complex-valued). Since then the posterior distribution $p(\mathbf{u}|\mathbf{v})$ is Gaussian, the posterior mean $\boldsymbol{\mu}_{a}$ equals the value that minimizes the quadratic $J(\mathbf{u})$ in \eqref{Posterior_Supp_SI}, which is given by,
\begin{equation}\label{Posterior_Mean_SI}
  \boldsymbol{\mu}_{a} = (\mathbf{I}-\mathbf{K}\mathbf{G})\boldsymbol{\mu}_f + \mathbf{K}\mathbf{v},
\end{equation}
where,
\begin{equation}\label{K_SI}
  \mathbf{K}= \mathbf{R}_{f}\mathbf{G}^\mathtt{T}(\mathbf{G}\mathbf{R}_{f}\mathbf{G}^\mathtt{T} + \mathbf{R}^o)^{-1},
\end{equation}
is known as the Kalman gain. The formulae in \eqref{Posterior_Mean_SI}--\eqref{K_SI} can be obtained by taking the derivative of $J(\mathbf{u})$ with respect to $\mathbf{u}$ using vector calculus and finding the root or solution of the resulting linear equation. Correspondingly, the posterior variance $\mathbf{R}_a$ is given by,
\begin{equation}\label{Posterior_Var_SI}
\mathbf{R}_a = (\mathbf{I}-\mathbf{K}\mathbf{G})\mathbf{R}_f.
\end{equation}

\paragraph{Posterior distribution as a function of the number of observations $L$.} Let us consider the scalar case of $u$, namely $u$ being a one-dimensional state variable. This is purely for the sake of simplicity, with the general case of $u$ being multi-dimensional and easily adaptable, mutatis mutandis. Despite $u$ being a scalar, $\mathbf{v}$ (and correspondingly $\boldsymbol{\epsilon}$) can be a vector, representing multiple measurements of $u$. In such a case, the observation $\mathbf{v}$ is an $L\times 1$ vector, where $L$ denotes the number of said measurements. For further notational simplicity, let us assume that the prior variance is $R_f=1$ and that the observational noise is $\mathbf{R}^o=\mathbf{I}$, which is an $L\times L$ identity matrix, meaning that the noise of each measurement is of unit variance and that the measurements are independent of each other (due to their normality). The observational operator $\mathbf{g}$ is an $L\times 1$ column vector. The dimension of the outer product $\mathbf{g}\mathbf{g}^\mathtt{T}$, which will be used below, is $L\times L$. With these, the Kalman gain $\mathbf{K}$ is necessarily a $1\times L$ row vector, since in \eqref{K_SI} the dimensions of $\mathbf{g}^\mathtt{T}$ and $(\mathbf{g}\mathbf{g}^\mathtt{T} + \mathbf{I})^{-1}$ are $1\times L$ and $L\times L$, respectively, and so we denote it by $\mathbf{k}^\mathtt{T}$ for the remainder of this subsection, and its expression can be simplified to,
\begin{equation}\label{K_Simplified_SI}
  \mathbf{k}^\mathtt{T} = \mathbf{g}^\mathtt{T}(\mathbf{g}\mathbf{g}^\mathtt{T} + \mathbf{I})^{-1}.
\end{equation}
In light of \eqref{K_Simplified_SI}, the corresponding expressions for the posterior mean \eqref{Posterior_Mean_SI} and posterior variance \eqref{Posterior_Var_SI} are,
\begin{equation}\label{Posterior_Mean_Simplified_SI}
  \mu_{a} = \frac{1}{\mathbf{g}^\mathtt{T}\mathbf{g} + 1}\mu_f + \frac{\mathbf{g}^\mathtt{T}}{\mathbf{g}^\mathtt{T}\mathbf{g} + 1} \mathbf{v},
\end{equation}
and,
\begin{equation}\label{Posterior_Var_Simplified_SI}
R_a = \frac{1}{\mathbf{g}^\mathtt{T}\mathbf{g} + 1}.
\end{equation}
Note that the $L\times L$ matrix $\mathbf{g}\mathbf{g}^\mathtt{T}$ appears in the parentheses in \eqref{K_Simplified_SI} while the denominators in \eqref{Posterior_Mean_Simplified_SI} and \eqref{Posterior_Var_Simplified_SI} contain the scalar inner product $\mathbf{g}^\mathtt{T}\mathbf{g}$, which is precisely the squared Euclidean norm of $\mathbf{g}$. Both $\mu_{a}$ and $R_a$ are scalars. The first term on the right-hand side of \eqref{Posterior_Mean_Simplified_SI} is a product of two scalars, while the second term is a product of a $1\times L$ row vector and an $L\times 1$ column vector.

The results in \eqref{Posterior_Mean_Simplified_SI}--\eqref{Posterior_Var_Simplified_SI} can be derived using the following two matrix identities (given by the well-known Woodbury matrix identity\footnote{Higham, Nicholas J. Accuracy and stability of numerical algorithms. Society for industrial and applied mathematics, 2002.}):
\begin{align}
(\mathbf{A}+\mathbf{B}\mathbf{C}\mathbf{D})^{-1} = \mathbf{A}^{-1}-\mathbf{A}^{-1}\mathbf{B}(\mathbf{C}^{-1} + \mathbf{D}\mathbf{A}^{-1}\mathbf{B})^{-1}\mathbf{D}\mathbf{A}^{-1},\label{Matrix_Identities_1_SI}\\
 (\mathbf{I}-(\mathbf{A}+\mathbf{E})^{-1}\mathbf{A})\mathbf{E}^{-1}=(\mathbf{A}+\mathbf{E})^{-1},\label{Matrix_Identities_2_SI}
\end{align}
where $\mathbf{A}$, $\mathbf{B}$, $\mathbf{C}$, $\mathbf{D}$, and $\mathbf{E}$ are matrices of sizes $m\times m$, $m\times L$, $L\times L$, $L\times m$, and $m\times m$, respectively.
According to \eqref{K_Simplified_SI},
\begin{equation}\label{1_KG_Simplified_SI}
  1-\mathbf{k}^\mathtt{T}\mathbf{g} = 1-\mathbf{g}^\mathtt{T}(\mathbf{g}\mathbf{g}^\mathtt{T} + \mathbf{I})^{-1}\mathbf{g}
\end{equation}
Denote by $\mathbf{A}=1$, $\mathbf{B}=\mathbf{g}^\mathtt{T}$, $\mathbf{C}=\mathbf{I}$ and $\mathbf{D}=\mathbf{g}$ (which essentially sets $m=1$ in the above matrix dimensions). The relationship in \eqref{Matrix_Identities_1_SI} simplifies \eqref{1_KG_Simplified_SI},
\begin{equation}\label{1_KG_Final_SI}
  1-\mathbf{k}^\mathtt{T}\mathbf{g} = \frac{1}{1+\mathbf{g}^\mathtt{T}\mathbf{g}},
\end{equation}
which is a scalar.
Likewise, denote by $\mathbf{A}=\mathbf{g}\mathbf{g}^\mathtt{T}$ and $\mathbf{E}=\mathbf{I}$, both of which are of size $L\times L$. The relationship in \eqref{Matrix_Identities_2_SI} leads to the $1\times L$ vector,
\begin{equation}\label{K_Final_SI}
\begin{split}
  \mathbf{k}^\mathtt{T}&= \mathbf{g}^\mathtt{T}(\mathbf{g}\mathbf{g}^\mathtt{T} + \mathbf{I})^{-1} = \mathbf{g}^\mathtt{T} (\mathbf{A}+\mathbf{E})^{-1} \\
  &= \mathbf{g}^\mathtt{T} (\mathbf{I}-(\mathbf{A}+\mathbf{E})^{-1}\mathbf{A})\mathbf{E}^{-1} \\
  &= \mathbf{g}^\mathtt{T}\big(\mathbf{I} - (\mathbf{g}\mathbf{g}^\mathtt{T} + \mathbf{I})^{-1}\mathbf{g}\mathbf{g}^\mathtt{T}\big)\\
  & = \big(1 - \mathbf{g}^\mathtt{T}(\mathbf{g}\mathbf{g}^\mathtt{T} + \mathbf{I})^{-1}\mathbf{g}\big)\mathbf{g}^\mathtt{T}\\
  & = \frac{\mathbf{g}^\mathtt{T}}{1+\mathbf{g}^\mathtt{T}\mathbf{g}},
\end{split}
\end{equation}
where the last equality is due to \eqref{1_KG_Simplified_SI} and \eqref{1_KG_Final_SI}. Therefore, \eqref{1_KG_Final_SI} and \eqref{K_Final_SI} give the posterior mean expression \eqref{Posterior_Mean_Simplified_SI} while \eqref{1_KG_Final_SI} leads to the expression of the posterior variance \eqref{Posterior_Var_Simplified_SI}.

To provide a more intuitive explanation of the results, further assume the observational operator $\mathbf{g}$ is given by an $L\times1$ vector $\mathbf{g} = (1,\ldots,1)^\mathtt{T}$, which means there are $L$ observations,
\begin{equation}\label{L_obs_SI}
  \mathbf{v} = \mathbf{g}u+\boldsymbol{\epsilon} \ \Leftrightarrow \ \left(
    \begin{array}{c}
      v_{1} \\
      \vdots \\
      v_{L} \\
    \end{array}
  \right) = \left(
    \begin{array}{c}
      1 \\
      \vdots \\
      1 \\
    \end{array}
  \right)u + \left(
    \begin{array}{c}
      \epsilon_{1} \\
      \vdots \\
      \epsilon_{L} \\
    \end{array}
  \right)
\end{equation}
where $\epsilon_1,\ldots\epsilon_L$ are $L$ independent standard Gaussian white noises, each having a variance of $1$. Thus, \eqref{Posterior_Mean_Simplified_SI} and \eqref{Posterior_Var_Simplified_SI} can be further simplified to
\begin{equation}\label{Posterior_Mean_Simplified_Final_SI}
  \mu_{a} = \frac{1}{L + 1}\mu_f + \frac{1}{L + 1}\sum_{l=1}^L v_l,
\end{equation}
and
\begin{equation}\label{Posterior_Var_Simplified_Final_SI}
R_a = \frac{1}{L + 1}.
\end{equation}

The result in \eqref{Posterior_Mean_Simplified_Final_SI} indicates that the observational noise can cancel off when multiple observations appear. This is consistent with the law of large numbers. In particular, when more and more observations are available, which means the noise is canceled off to a more significant degree, the weight toward the prior information becomes less and less. Another way to express this observation is to notice that the posterior mean in \eqref{Posterior_Mean_Simplified_SI} can be expressed as a convex linear combination of the prior mean and the maximum likelihood estimator (MLE). Specifically, observe that for,
\begin{equation} \label{Weight_Combo}
    w = \frac{1}{\mathbf{g}^\mathtt{T}\mathbf{g}+1},
\end{equation}
we have from \eqref{Posterior_Mean_Simplified_SI} and \eqref{L_obs_SI} that,
\begin{equation}\label{Posterior_Mean_Convex_Combo}
    \mu_a = \frac{1}{\mathbf{g}^\mathtt{T}\mathbf{g} + 1}\mu_f + \frac{\mathbf{g}^\mathtt{T}}{\mathbf{g}^\mathtt{T}\mathbf{g} + 1} \mathbf{v}  = w\mu_f + (1-w) \frac{\mathbf{g}^\mathtt{T}\mathbf{v}}{\mathbf{g}^\mathtt{T}\mathbf{g}} = \frac{1}{L + 1}\mu_f + \frac{L}{L + 1}\frac{\sum_{l=1}^L v_l}{L}.
\end{equation}
This is a general property of exponential families of distributions, where when assuming an exponential family parameterized by its natural parameter (i.e.\ in canonical form) and a conjugate family for a prior (necessarily an exponential family of distribution as well), then the posterior expectation of the mean of the natural statistic is a convex combination of the prior expectation of said natural statistic and the MLE\footnote{Diaconis, Persi and Ylvisaker, Donald. Conjugate Priors for Exponential Families. The Annals of Statistics, Ann. Statist. 7(2), 269-281, 1979.}. Once again, \eqref{Posterior_Mean_Convex_Combo} showcases the observation that as the number of observations increases, or $L\to\infty$, then more trust is put onto the MLE instead of our prior beliefs, i.e., more weight is put onto the observations which pull the state estimation towards them, and hopefully towards the underlying truth.

\begin{figure}[ht!]
	\begin{center}
		\hspace*{-0.5cm}\includegraphics[width=17.5cm]{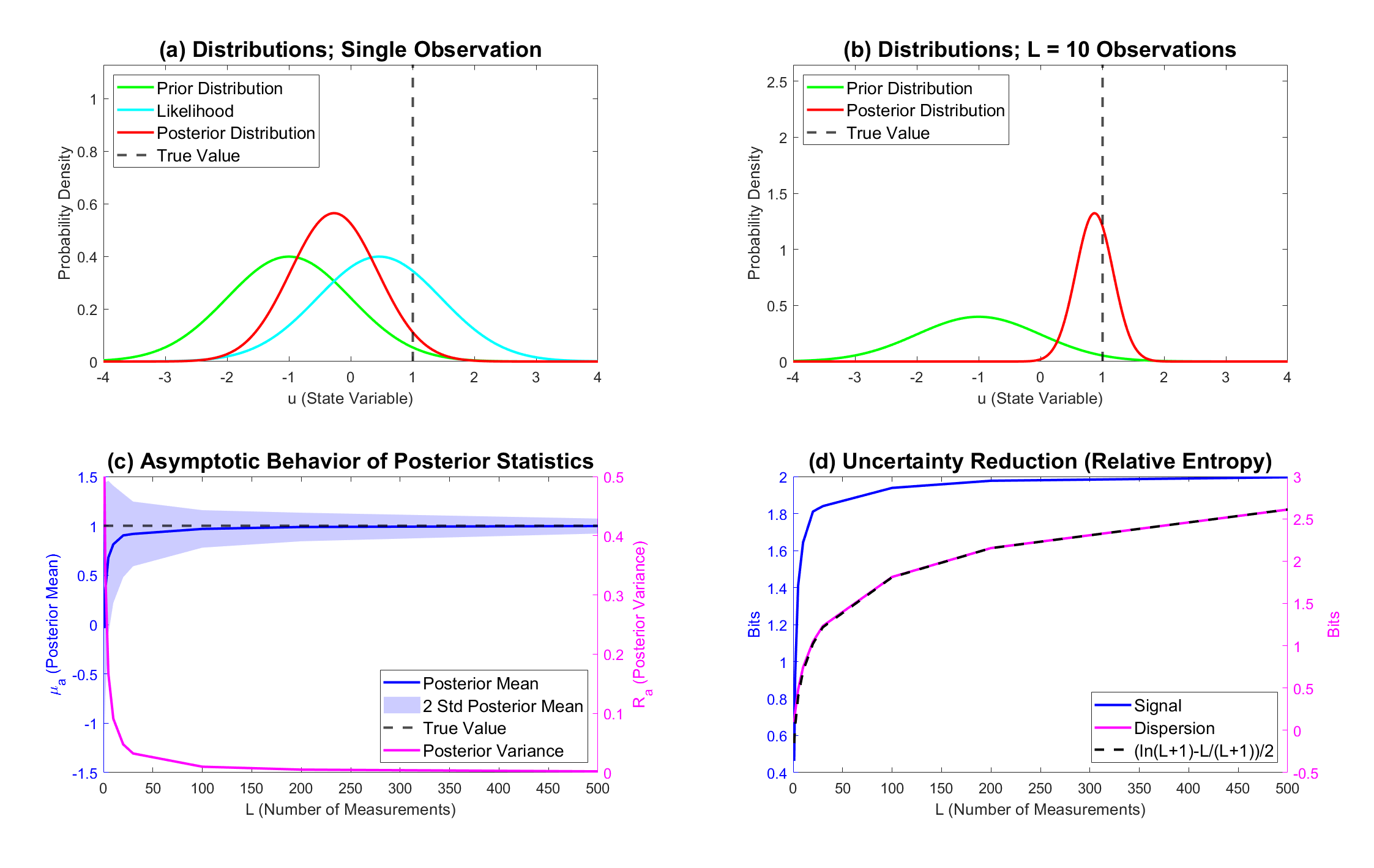}
	\end{center}
	\caption{Numerical illustration of the Bayes formula. Here $\mu_f=1$, $R_f=1$, and $m=1$. The observational noise is assumed to be $1$ for each individual observation. Panels (a)--(b): PDFs with different number of observations $L$. Panel (c): asymptotic behavior of the posterior mean $\mu_a$ and the posterior variance $R_a$ as a function of $L$. Due to the randomness in observations, the shaded area shows the variation of the results with 100 sets of independent observations for each fixed $L$. Panel (d): signal and dispersion components of the relative entropy, based on the solid curves in Panel (c). The dashed curve shows the theoretic value given by \eqref{Dispersion_Multiple_Obs_SI}, which is consistent with the numerical simulation.   }
	\label{SI_Bayes}
\end{figure}

These multiple observations also help reduce the uncertainty reflected by the posterior variance. One interesting fact is the following. By plugging \eqref{Posterior_Var_Simplified_Final_SI} into the dispersion part of the relative entropy expression \eqref{Relative_Entropy_Gaussian_SI}, it is shown that,
\begin{align}
\begin{split}
    \mbox{Dispersion~}&=\frac{1}{2}\Big[R_a(R_f)^{-1}-1-\ln(R_a(R_f)^{-1})\Big]\\
    &= \frac{1}{2}\left(R_a-1-\ln R_a\right) = \frac{1}{2}\left(-\frac{L}{L+1} + \ln(L+1)\right),
\end{split} \label{Dispersion_Multiple_Obs_SI}
\end{align}
which asymptotically, as $L\to\infty$, behaves like $\ln L$. Despite being extremely simple, this example provides an intuition of the UQ in the uncertainty reduction of the LaDA in the main text, explaining why the information reduction is a function of the logarithm of the tracer number. Figure \ref{SI_Bayes} numerically confirms this finding.

\subsection{Lagrangian DA}
{\color{blue}{
The codes that correspond to the content of this subsection are:
\begin{itemize}
  \item MATLAB codes: \textbf{Flow\_Model.m}, \textbf{LDA\_Main\_Filter.m}, and \textbf{LDA\_Function\_of\_L.m}. The code \textbf{Flow\_Model.m} needs to be run first to generate the true flow field. Then \textbf{LDA\_Function\_of\_L.m} should be run to obtain the posterior distribution and the uncertainty reduction as a function of $L$. It will call the script \textbf{LDA\_Main\_Filter.m} automatically.
  \item Python codes: \textbf{Flow\_Model.py}, \textbf{LDA\_Main\_Filter.py}, and \textbf{LDA\_Function\_of\_L.py}. The code \textbf{Flow\_Model.py} needs to be run first to generate the true flow field. Then \textbf{LDA\_Function\_of\_L.py} should be run to obtain the posterior distribution and the uncertainty reduction as a function of $L$. It will call the script \textbf{LDA\_Main\_Filter.py} automatically.
\end{itemize}}}

\noindent The LaDA scheme consists of two sets of equations: One for the observational process and the other for the flow field model.

\paragraph{Forecast Model.}  Consider a random flow modeled by a finite number of Fourier modes with random amplitudes in the periodic domain $(-\pi,\pi]^2$,
\begin{equation}\label{LaDA_Velocity_Field_SI}
    \mathbf{u}({\mathbf{x}},t) = \sum_{\mathbf{k}\in \mathbf{K}} {\hat{u}_{\mathbf{k}}(t)}\cdot e^{i\mathbf{k}\cdot {\mathbf{x}}}\cdot\mathbf{r}_{\mathbf{k}},
\end{equation}
where $\mathbf{r}_{\mathbf{k}}$ is the so-called eigenvector that contains the relationship between different components in $\mathbf{u}$. For example, if the flow field $\mathbf{u}$ is incompressible, then the divergence-free condition is reflected in the eigenvector of each Fourier mode.
As a forecast model, let us adopt a linear stochastic model (known as a complex-valued Ornstein--Uhlenbeck (OU) process)
as a surrogate to describe the time evolution of each ${\hat{u}_{\mathbf{k}}(t)}$,
\begin{equation}\label{LaDA_LSM_SI}
    \frac{\d{\hat{u}_{\mathbf{k}}(t)}}{\d t}=(-d_{\mathbf{k}} + i\omega_{\mathbf{k}}) {\hat{u}_{\mathbf{k}}(t)}+f_{\mathbf{k}}(t)+{ \sigma_{\mathbf{k}}\dot{W}^u_{\mathbf{k}}(t)}.
\end{equation}
Note that the governing equation of each Fourier coefficient $\hat{u}_\mathbf{k}(t)$, as was discussed in the ``UQ in Advancing Efficient Modeling'' of the main text, usually has the following nonlinear form,
\begin{equation}\label{LaDA_Original_NonlinearEq_SI}
    \frac{\d{\hat{u}_\mathbf{k}(t)}}{\d t} = \big(-d_\mathbf{k} +i\omega_\mathbf{k}\big){\hat{u}_\mathbf{k}(t)}+{\hat{f}_\mathbf{k}(t)}+{\sum_{\mathbf{m}\in\mathbf{K}} c_{\mathbf{m},\mathbf{k}}\hat{u}_\mathbf{m}\hat{u}_{\mathbf{k}-\mathbf{m}}},
\end{equation}
where the convolution term is due to the nonlinearity in the original PDE model when passed through a spectral discretization. Since many nonlinear terms in the original equation involve high frequencies (e.g., $\hat{u}_{\mathbf{m}}$ with $|\mathbf{m}|\gg |\mathbf{k}|$), the time series constructed by these terms is fully chaotic and has very short temporal correlations. This justifies replacing the complicated nonlinearity with effective stochastic noise (and additional damping) in the linear stochastic process as a statistical forecast model. It is consistent with the stochastic mode reduction --- the equations of motion for the unresolved fast modes are modified by representing the nonlinear self-interaction terms between unresolved modes by damping and stochastic terms.

\paragraph{Observations.} The observations are given by the trajectories of $L$ tracers. Each of the tracer displacements $\mathbf{x}_l=(x_l,y_l)^\mathtt{T}$ satisfy a stochastic perturbation of the usual velocity relation as to account for possible measurement errors and the contribution from the velocity from the unresolved scales (with constant noise under the assumption of the tracers being modelled as point particles with no mass). That is,
  \begin{equation}\label{LaDA_Observations_SI}
    \frac{\d{\mathbf{x}_l(t)}}{\d t}={\mathbf{u}}({\mathbf{x}_l(t)}, t)+\sigma_x \dot{\mathbf{W}}^x_{l}(t)
    =\sum_{\mathbf{k}\in \mathbf{K}} \underbrace{{\hat{u}_{\mathbf{k}}(t)} \cdot e^{i\mathbf{k}\cdot {\mathbf{x}_l(t)}}\cdot\mathbf{r}_{\mathbf{k}}}_{\mbox{Nonlinear!}}+\sigma_x \dot{\mathbf{W}}^x_l(t),\quad l=1,\ldots, L.
\end{equation}
Note that the observational process is a highly nonlinear equation of the variable $\mathbf{x}$ due to its appearance in the exponential function on the right-hand side. Furthermore, since the tracers are treated as massless particles, the velocity of each one coincides with that of the underlying velocity field at that location.

\paragraph{Coupled system for the LaDA.}
Define
\begin{equation}\label{LaDA_Abstract_Variables_SI}
{\hat{\mathbf{u}}=\big(\hat{u}_{\mathbf{k}_1},...\hat{u}_{\mathbf{k}_{|\mathbf{K}|}}\big)^{\mathtt{T}}},\qquad\mbox{and}\qquad{\mathbf{x}=\big(x_{1},y_{1},...,x_{L},y_{L}\big)^{\mathtt{T}}}.
\end{equation}
The abstract form of the above coupled system \eqref{LaDA_Velocity_Field_SI}--\eqref{LaDA_Observations_SI} for the LaDA is given by
\begin{equation}\label{LaDA_Coupled_System_SI}
\begin{split}
    \mbox{Observational processes:}\qquad \frac{\d{ \mathbf{x}}}{\d t} &= \mathbf{P}({\mathbf{x}}){\hat{\mathbf{u}}}+\boldsymbol{\Sigma}_{\mathbf{x}} \dot{\mathbf{W}}_{\mathbf{x}},\\
    \mbox{Forecast flow model:}\qquad \frac{\d\hat{\mathbf{u}}}{\d t} &= -\boldsymbol\Gamma {\hat{\mathbf{u}}}+\mathbf{f}+\boldsymbol{\Sigma}_{\hat{\mathbf{u}}} \dot{\mathbf{W}}_{\hat{\mathbf{u}}}.
\end{split}
\end{equation}
The LaDA (filtering) solution, namely the posterior distribution, is given by $p({ \hat{\mathbf{u}}(t)}|{ \mathbf{x}(s\leq t)})$, and is optimal in the sense that it minimizes the mean squared error. This is a nonlinear DA problem since the coupled system for $({\mathbf{x}},{\hat{\mathbf{u}}})$ is nonlinear. Nevertheless, conditioned on the observed tracer trajectories ${\mathbf{x}(s\leq t)}$, the flow variable ${\hat{\mathbf{u}}(t)}$ appears linearly in both the observational processes and the forecast model. Therefore, the posterior distribution $p({ \hat{\mathbf{u}}(t)}|{ \mathbf{x}(s\leq t)})\sim\mathcal{N}(\boldsymbol\mu,\mathbf{R})$ is Gaussian. Systems of this form, for which the conditional distribution of the unobservable or latent process is Gaussian when conditioning on the partial observations, are known as conditional Gaussian nonlinear systems \cite{chen2023stochastic}. Notably, despite the nonlinearity, closed analytic formulae are available for computing $p({ \hat{\mathbf{u}}(t)}|{ \mathbf{x}(s\leq t)})$, which allows an exact and accurate LaDA solution \cite{chen2023stochastic},
\begin{equation}\label{LaDA_Filter_SI}
\begin{split}
\frac{\d\boldsymbol{\mu}}{\d t} &= \left(\mathbf{f} - \boldsymbol{\Gamma} \boldsymbol{\mu}\right)  + \mathbf{R}\mathbf{P}^\ast(\boldsymbol\Sigma_{\mathbf{x}}\boldsymbol\Sigma_{\mathbf{x}}^*)^{-1}\left(\frac{\d \mathbf{x}}{\d t} - \mathbf{P}\boldsymbol{\mu} \right),\\
\frac{\d\mathbf{R}}{\d t} &= -\boldsymbol{\Gamma}\mathbf{R} - \mathbf{R}\boldsymbol{\Gamma}^\ast + \boldsymbol{\Sigma}_{\hat{\mathbf{u}}}\boldsymbol{\Sigma}_{\hat{\mathbf{u}}}^\ast - \mathbf{R}\mathbf{P}^\ast(\boldsymbol\Sigma_{\mathbf{x}}\boldsymbol\Sigma_{\mathbf{x}}^*)^{-1}\mathbf{P}\mathbf{R}.
\end{split}
\end{equation}

\paragraph{Numerical illustration.} In the numerical tests, the true signal is assumed to be generated from the linear stochastic model \eqref{LaDA_LSM_SI} for simplicity. In practice, the true signal can be generated from other models or direction observations. The underlying dynamics of the true signal do not need to be known. Only the resulting Lagrangian trajectories are required as input for the LaDA.

\begin{figure}[ht!]
	\begin{center}
		\hspace*{-0.5cm}\includegraphics[width=17cm]{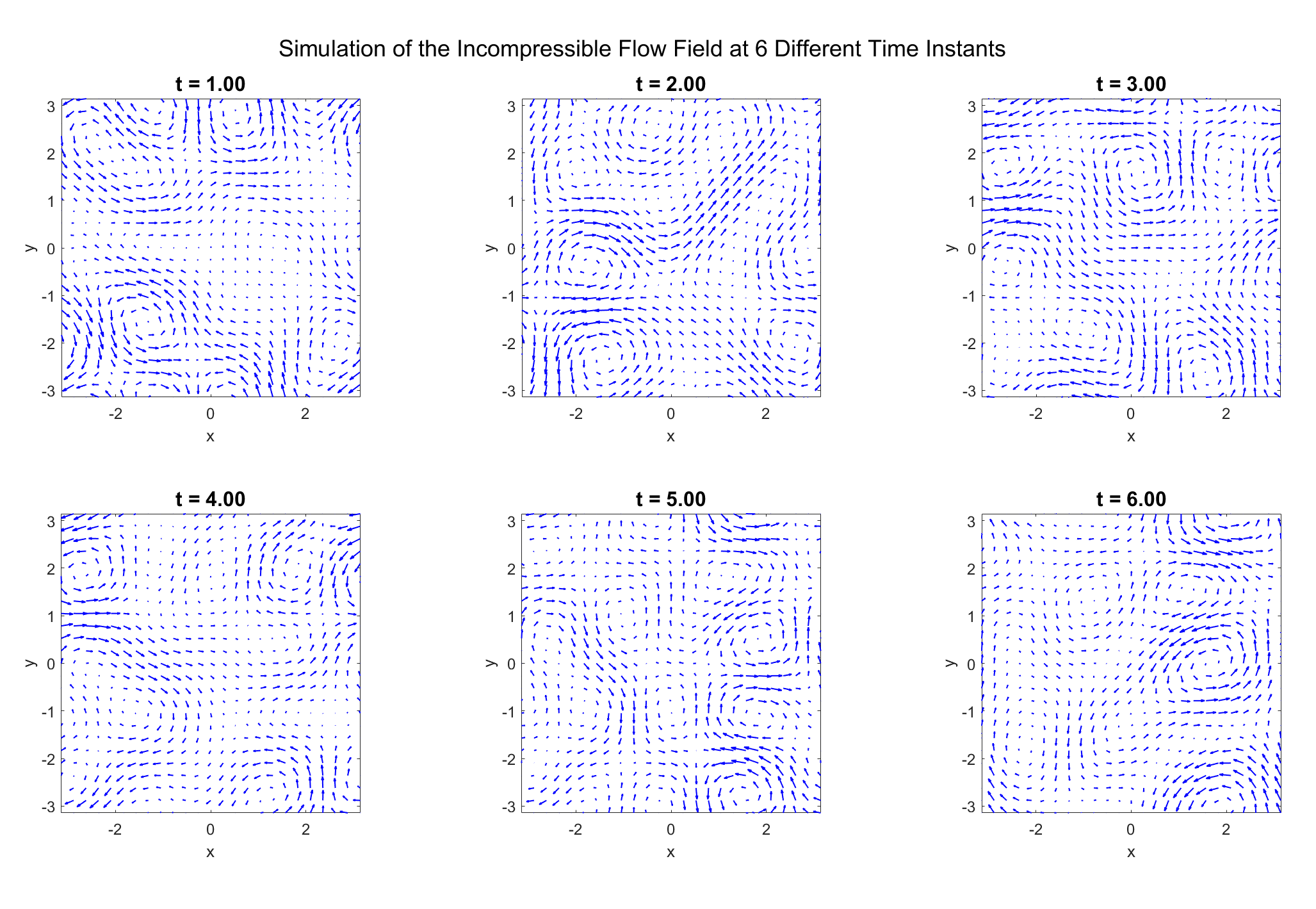}
	\end{center}
	\caption{Simulation of the flow field at different times. These flow fields will be used to generate the observed Lagrangian tracers.   }
	\label{SI_LaDA_Flow_Field}
\end{figure}

In the simulation here, the Fourier wavenumbers are $\mathbf{k}\in[-2,2]^2\cap \mathbb{Z}^2$, excluding the mode in the origin $(0,0)$ which is the mode that corresponds to a mean sweep or background flow. The parameters for all these 24 modes are the same, and are given by $d_\mathbf{k}=0.5$, $\omega_\mathbf{k}=0$, $f_\mathbf{k}$ and $\sigma_\mathbf{k}=0.5$ for all $\mathbf{k}$ (with these values the coefficients satisfy the usual conjugate or reality condition needed between conjugate wavenumbers as to preserve the real values of the underlying flow field). Figure \ref{SI_LaDA_Flow_Field} shows the snapshots of the true flow field, which is one random realization from the flow model. The observed Lagrangian tracer trajectories are driven by this flow field.

\begin{figure}[ht!]
	\begin{center}
		\hspace*{-0.5cm}\includegraphics[width=17.5cm]{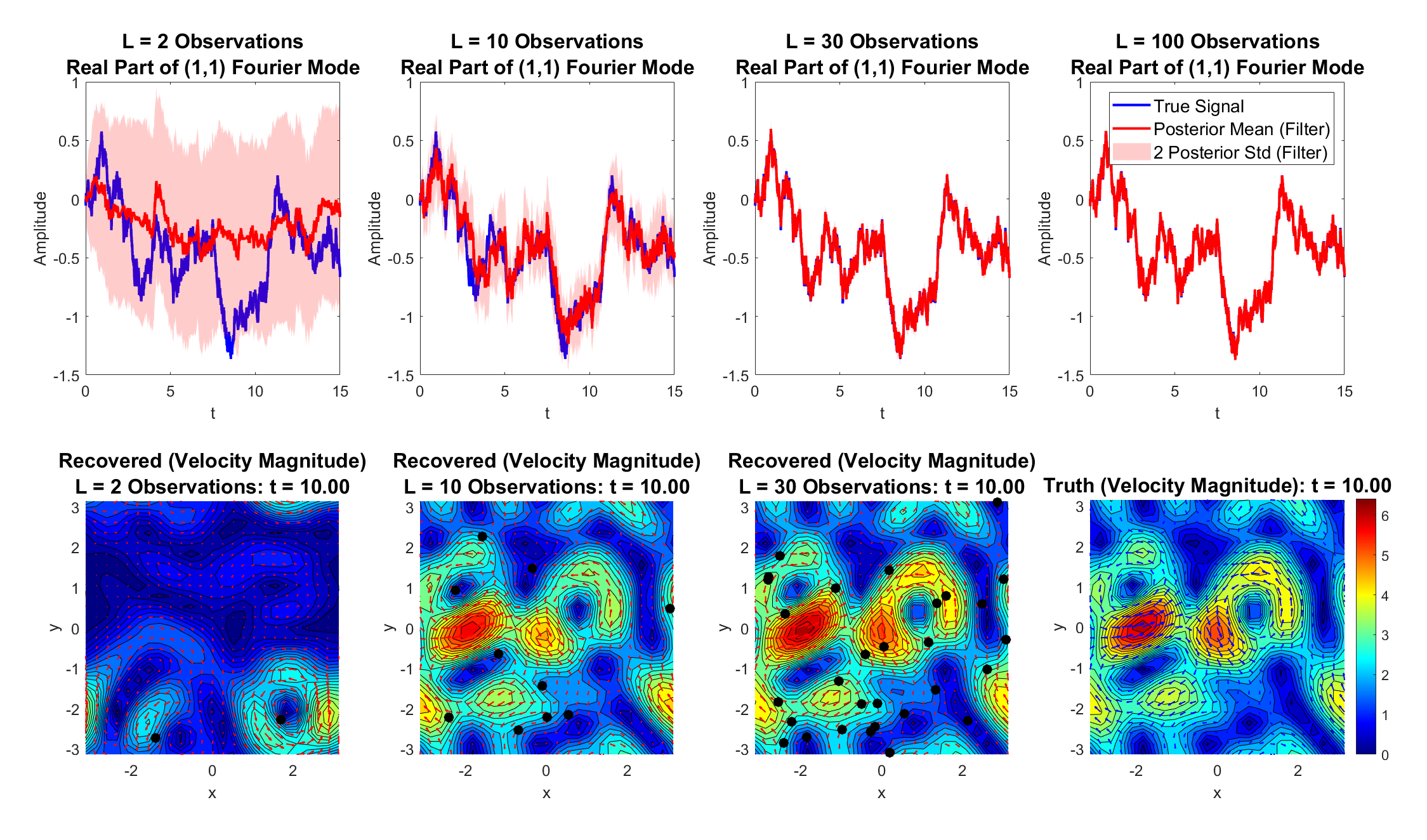}
	\end{center}
	\caption{Recovering the flow fields from the LaDA. Top row: The recovery of one Fourier mode $(1,1)$ (red) with different observed tracers $L$ compared with the truth (blue). The `std' stands for the standard deviation. Bottom row: The recovered flow field reconstructed from the recovered Fourier modes (posterior mean) compared with the truth (last panel). }
	\label{SI_LaDA_Recovered_with_L}
\end{figure}

Figure \ref{SI_LaDA_Recovered_with_L} shows the results from the LaDA. The top row shows the recovery of one Fourier mode $(1,1)$ (red) with different observed tracers $L$ compared with the truth (blue). It is seen that when $L$ becomes large, the posterior mean converges to the truth, and the posterior uncertainty shrinks, consistent with the results shown in the main text. The bottom row illustrates the recovered flow field reconstructed from the recovered Fourier modes (posterior mean) and is compared with the truth (last panel), which validates the convergence of the posterior mean estimate.

\section{Role of the Uncertainty in Diagnostics}
\subsection{Parameter estimation with uncertainties in data}

{\color{blue}{
The codes that correspond to the content of this subsection are:
\begin{itemize}
  \item MATLAB code: \textbf{Parameter\_Estimation.m}
  \item Python code: \textbf{Parameter\_Estimation.py}
\end{itemize}}}

\noindent In this section, we will use numerical simulations to illustrate more intuitively how the uncertainty affects parameter estimation.

Consider the same linear system as in the main text,
\begin{equation}\label{Linear_Oscillator_SI}
  \frac{\d x}{\d t} = ay,\qquad \frac{\d y}{\d t} = bx,
\end{equation}
Define a matrix
\begin{equation}\label{Matrix_M_SI}
\mathbf{M}_i=\begin{pmatrix}
y_i & 0 \\
0 & x_i \\
\end{pmatrix},
\end{equation}
which takes values at time $t_i$. As was discussed in the main text, if the observational data of $x$ and $y$ and their time derivatives $\dot{x}:=\d x/\d t$ and $\dot{y}:=\d y/\d t$ are available at time $t_i$ for $i=1,\ldots, I$, then the parameter $a$ can be estimated via standard linear regression and through the least-squares solution
\begin{equation}\label{regression_parameter_estimation_a_SI}
  a = \left(\sum_i y_i^2\right)^{-1}\left(\sum_i y_i\dot{x}_i\right).
\end{equation}
However, when only $x_i$ and $\dot{x}_i$ are observed, the expectation needs to be taken to consider the uncertainty in the unobservable, $y$. Consequently, the parameter estimation formula needs to be modified to,
\begin{equation}\label{regression_parameter_estimation_a_UQ_SI}
  a = \left(\sum_i \left(\langle y_i^2\rangle + \langle(y_i^\prime)^2\rangle\right)\right)^{-1}\left(\sum_i \langle y_i\rangle\dot{x}_i\right),
\end{equation}
where $y_i$ has been written into the Reynolds decomposition form as $y_i = \langle y_i\rangle + y_i^\prime$ and the term $\langle(y_i^\prime)^2\rangle$ in \eqref{regression_parameter_estimation_a_UQ_SI} stems from averaging the nonlinear function $y_i^2$ in \eqref{regression_parameter_estimation_a}.

\begin{figure}[ht!]
	\begin{center}
		\hspace*{-0.5cm}\includegraphics[width=17.5cm]{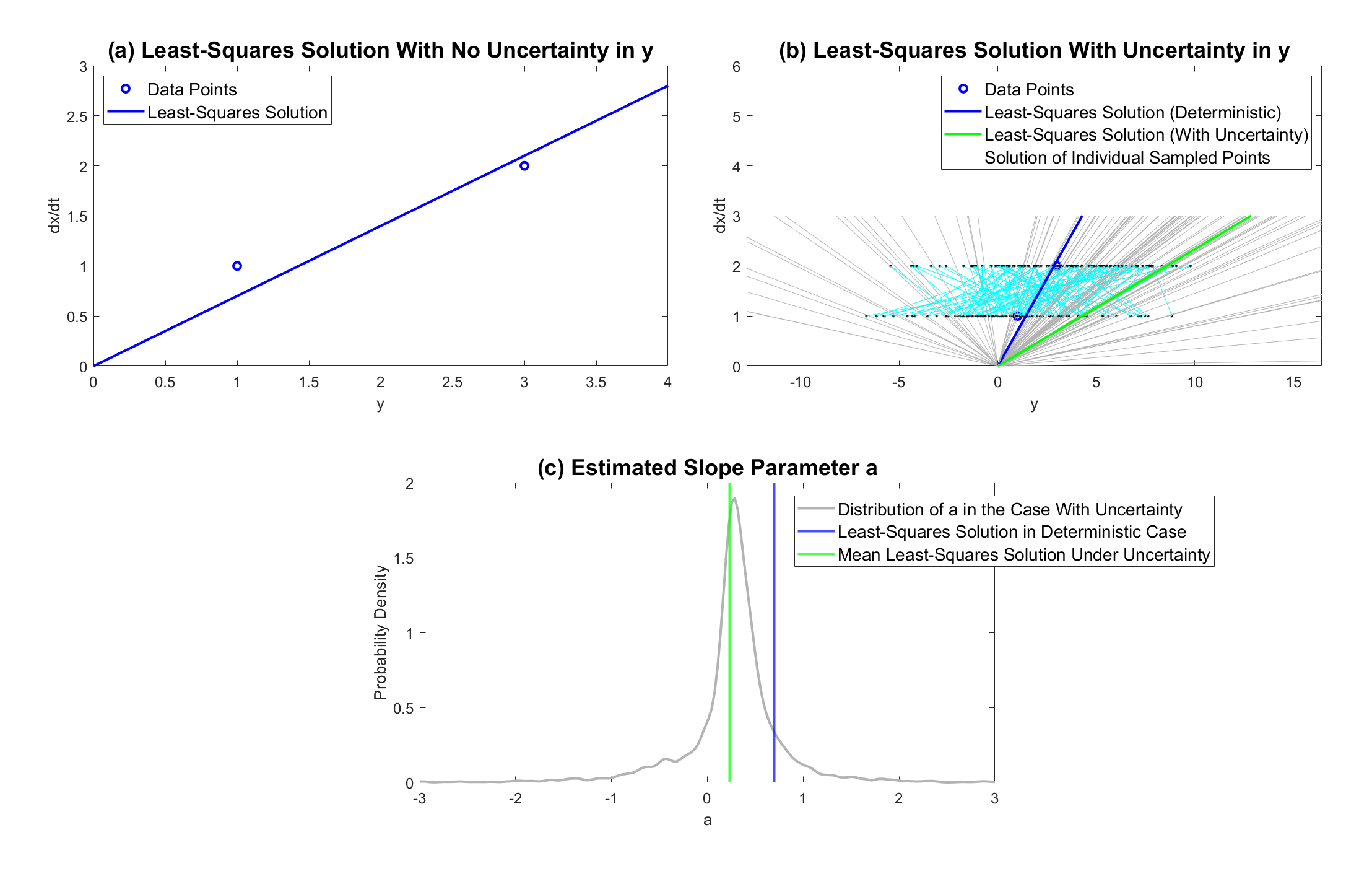}
	\end{center}
	\caption{Parameter estimation with two data points $(\dot{x}_1, y_1) = (1,1)$ and $(\dot{x}_2, y_2) = (2,3)$. Panel (a): When the variables $y_1$ and $y_2$ are observed, applying the linear regression \eqref{regression_parameter_estimation_a_SI} gives the least-squares solution. Panel (b): Each black dot is one sampled point, assuming $y_i$ ($i=1, 2$) satisfy a Gaussian distribution centered at the true value with a large uncertainty (variance) $r=10$. Each cyan line indicates the two points of $y_1$ and $y_2$ as a one sample pair. One gray line is the estimated parameter $a$ for one sample pair of $(y_1,y_2)$. The average value is shown in green, which is different from that in the deterministic case without uncertainty (blue). Panel (c): The distribution of $a$ and its average value, compared with the one in the deterministic case without uncertainty (blue). }
	\label{SI_Parameter_Estimation}
\end{figure}

\paragraph{Numerical validation.} Consider the following simple example with two data points: $(\dot{x}_1, y_1) = (1,1)$ and $(\dot{x}_2, y_2) = (2,3)$. Panel (a) in Figure \ref{SI_Parameter_Estimation} shows the standard linear regression solution for estimating the parameter $a$ in \eqref{Linear_Oscillator_SI} when the data of $\dot{x}_i$ and $y_i$ are both available.

Next, consider the case where $y_i$ is not directly observed. Assume $y_i$ satisfies a Gaussian distribution, characterizing its uncertainty. The Gaussian distribution for $y_1$ and $y_2$ is centered at the true value but has a large variance of $r = 10$. The formula \eqref{regression_parameter_estimation_a_UQ_SI} can be directly applied to estimate $a$ when such an uncertainty in $y_i$ appears. The results are shown in the green line in Panels (b)--(c), which is different from the blue line that corresponds to the solution of the deterministic case.

In addition to using \eqref{regression_parameter_estimation_a_UQ_SI}, we also show the results by sampling from the distribution of $y_i$, estimating $a$ based on each pair of samples, and then computing the average value of $a$. This is the motivation behind \eqref{regression_parameter_estimation_a_UQ_SI}. When the number of samples becomes infinity, the numerical results converge to the theoretical value in \eqref{regression_parameter_estimation_a_UQ_SI}.

In Panel (b), the black dots are the samples for $y_1$ and $y_2$, sampled independently. Each cyan line indicates the two points of $y_1$ and $y_2$ as a one sample pair. The true values of $y_1$ and $y_2$ with $y_2=y_1+2$ are seen by the blue dots. However, due to the large uncertainty, it is likely that $y_2<y_1$ or $y_2\gg y_1$ in many sample pairs. Therefore, the value $a$ from the least-squares solution, which is the slope of the gray lines, varies for different samples, where its distribution is shown in Panel (c). The average value of $a$, shown in green, differs from the value in the deterministic case without uncertainty (shown in blue).

As a final remark, it is worth noting that the distribution of $a$ shown in Panel (c) is highly non-Gaussian with fat tails. This is again counterintuitive since the underlying system is linear, and the uncertainty in $y_i$ is assumed to be Gaussian. Nevertheless, as the diagnostic \eqref{regression_parameter_estimation_a_UQ_SI} contains nonlinearity through $\langle (y_i^\prime)^2\rangle$, the non-Gaussian statistics naturally appear. In fact, \eqref{regression_parameter_estimation_a_UQ_SI} is effectively proportional to the reciprocal of $y_i$, which for a Gaussian random variable is heavy-tailed (and even bimodal). Therefore, as was discussed in the main text, the Gaussian distribution of $y_i$ after such a nonlinear transform becomes highly non-Gaussian.

\subsection{Eddy identification}

{\color{blue}{
The codes that correspond to the content of this subsection are:
\begin{itemize}
  \item MATLAB codes: \textbf{Flow\_Model.m}, \textbf{LDA\_Main\_Smoother.m}, and \textbf{Eddy\_Identification.m}. When running \textbf{Eddy\_Identification.m}, the other two scripts will be called automatically.
  \item Python codes: \textbf{Flow\_Model.py}, \textbf{LDA\_Main\_Smoother.py}, and \textbf{Eddy\_Identification.py}. When running \textbf{Eddy\_Identification.py}, the other two scripts will be called automatically.
\end{itemize}}}

\noindent In this article, the OW parameter is adopted as the criterion for identifying eddies. Recall that the OW parameter is defined as
\begin{equation}\label{OW_Parameter_SI}
        \operatorname{OW} = s_\mathrm{n}^2 + s_\mathrm{s}^2 - \omega^2,
    \end{equation}
    where the normal strain, the shear strain, and the relative vorticity are given by
    \begin{equation}
        s_\mathrm{n} = u_x-v_y, \quad s_\mathrm{s} = v_x+u_y, \quad\mbox{and}\quad\omega = v_x-u_y,
    \end{equation}
    respectively, with $\mathbf{u}=(u,v)$ being the two-dimensional velocity field and where the shorthand notation for the partial derivatives is being used $u_x:=\partial u/\partial x$. When the OW parameter is negative, the relative vorticity is larger than the strain components, indicating vortical flow. As was stated in the main text, the primary goal of this article is to show that nonlinearity in the eddy diagnostic will make UQ play a crucial role. Providing a rigorous definition of an eddy and discussing the pros and cons of different eddy identification methods are not the main focus here. A more detailed discussion of eddy definitions and the study of the uncertainty in affecting the eddy identification can be instead found in \cite{covington2024probabilistic}.

\begin{figure}[ht]
	\begin{center}
		\hspace*{-0.5cm}\includegraphics[width=17.5cm]{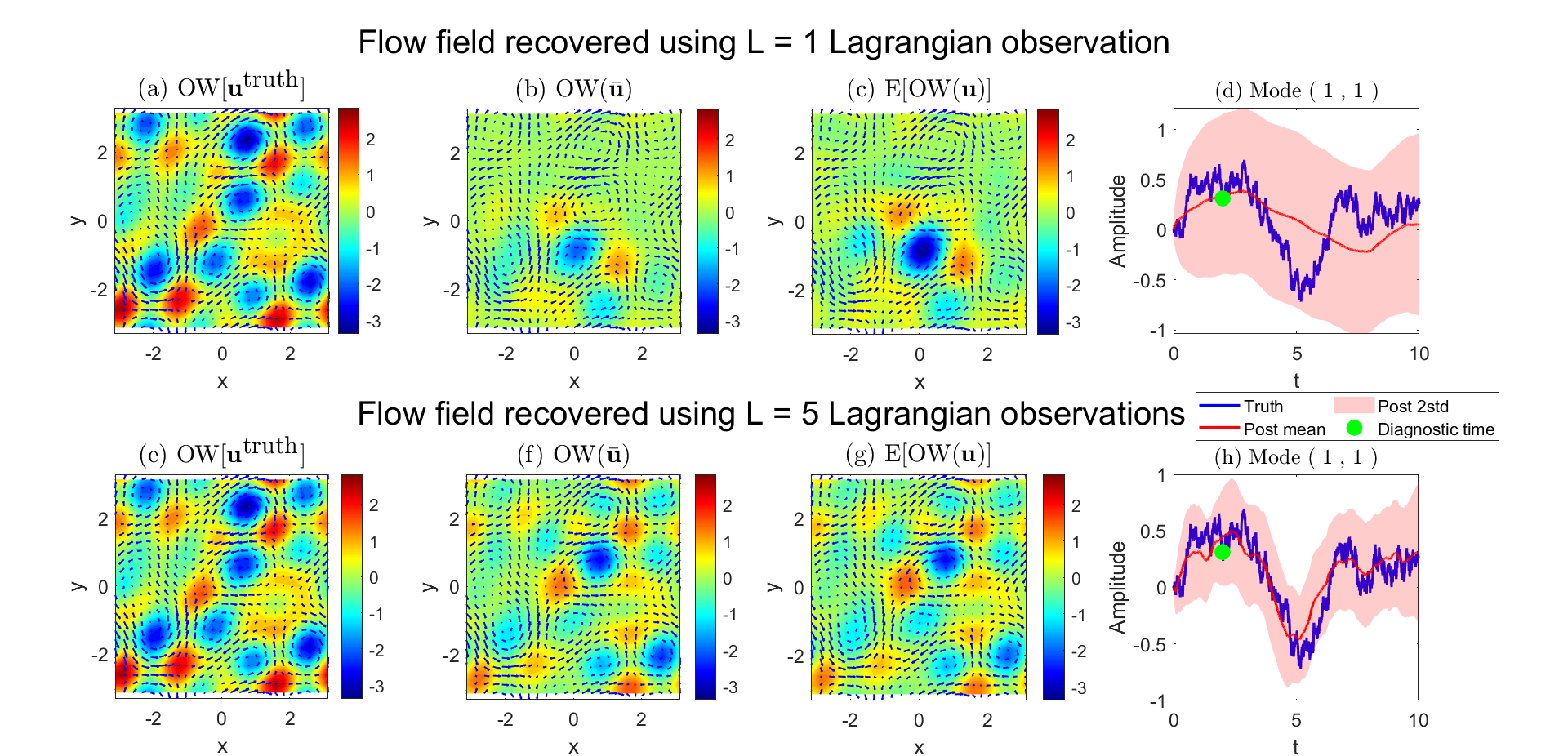}
	\end{center}
	\caption{Comparison of the OW parameters based on the recovered flow field using $L=1$ and $L=5$ observed Lagrangian tracers.  }
	\label{SI_Eddy_Comparison}
\end{figure}

Using the same setup of the flow field as in the LaDA section above (Section 4.2), Figure \ref{SI_Eddy_Comparison} compares the OW parameters based on the recovered flow field when using $L=1$ and $L=5$ observed Lagrangian tracers. When $L=1$, the uncertainty in the recovered flow field is large (Panel (d)). As a result, $\mbox{OW}(\bar{\mathbf{u}})$ (Panel (b)) and $\mbox{E[OW}({\mathbf{u}})]$ (Panel (c)) are different from the truth. Nevertheless, as was discussed in the main text, the nonlinearity in the OW parameter leads to the difference in $\mbox{OW}(\bar{\mathbf{u}})$ and $\mbox{E[OW}({\mathbf{u}})]$. Notably, both $\mbox{OW}(\bar{\mathbf{u}})$ and $\mbox{E[OW}({\mathbf{u}})]$ significantly miss the correct identification of many eddies due to the large uncertainty. In contrast, when the uncertainty is reduced by using more observations $L=5$, the results for $\mbox{OW}(\bar{\mathbf{u}})$ and $\mbox{E[OW}({\mathbf{u}})]$ are improved.

\begin{figure}[ht]
	\begin{center}
		\hspace*{-1.5cm}\includegraphics[width=20cm]{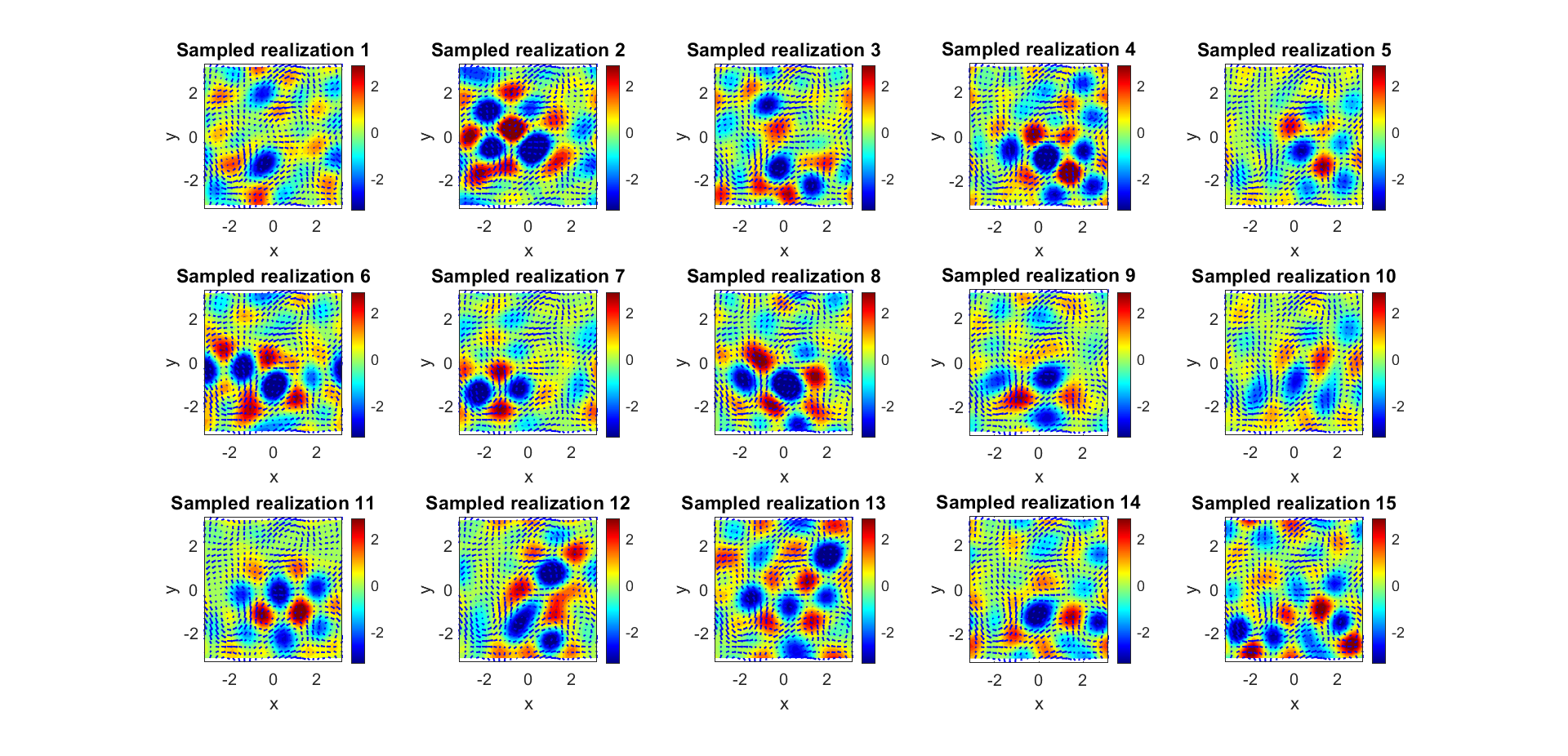}
	\end{center}
	\caption{Eddy identification based on the sampled realization of the flow field from the LaDA with $L=1$, corresponding to the case in the top row of Figure \ref{SI_Eddy_Comparison}. }
	\label{SI_Eddy_Sampling}
\end{figure}

Figure \ref{SI_Eddy_Sampling} shows the OW parameter based on the sampled realization of the flow field from the LaDA with $L=1$, corresponding to the case in the top row of Figure \ref{SI_Eddy_Comparison}. For different realizations, the resulting OW parameters are different. Since the uncertainty is large using when using only a single tracer observation ($L=1$), the difference in the resulting OW parameters is also very significant. Nevertheless, unlike $\mbox{OW}(\bar{\mathbf{u}})$ and $\mbox{E[OW}({\mathbf{u}})]$ that significantly underestimate the strength of the OW parameter, the OW parameters in all the panels in Figure \ref{SI_Eddy_Sampling} have comparable amplitude as the truth. Although these sampled realizations of the flow field are not the same as the truth, the resulting OW parameter at each location can be collected to construct a PDF describing the possible range of the OW parameter, indicating the probability of the occurrence of a specific eddy. Notably, these PDFs can also help calculate the distribution of the lifetime and size of each eddy. When uncertainty appears, such UQ in the diagnostics is important as the deterministic solution by averaging out the uncertainty in one way or the other (e.g., $\mbox{OW}(\bar{\mathbf{u}})$ and $\mbox{E[OW}({\mathbf{u}})]$) may lose a large amount of crucial information. See \cite{covington2024probabilistic} for details.

\section{Calibrating Stochastic Models Based on UQ}
{\color{blue}{
The codes that correspond to the content of this subsection are:
\begin{itemize}
  \item MATLAB code: \textbf{Calibrating\_Stochastic\_Model\_with\_UQ.m}
  \item Python code: \textbf{Calibrating\_Stochastic\_Model\_with\_UQ.py}
\end{itemize}}}

\noindent Consider a real-valued time series, possibly stemming from the realization of a dynamical stochastic system with highly nonlinear dynamics and/or multiplicative and additive noise. This section aims to (a) illustrate how to find the linear stochastic model that fits this time series by matching the key statistics and (b) the dynamical and statistical performance of such a calibrated linear stochastic model. The three statistics adopted for model calibration are the equilibrium mean, the equilibrium variance, and the decorrelation time of the given time series.

The real-valued scalar linear stochastic model that will be used for calibration is given by a simple OU process as,
\begin{equation}\label{LSM_SI}
  \frac{\d x}{\d t} = (- a x + f) \d t + \sigma\dot{W},
\end{equation}
where $\dot{W}$ is a white noise. The long-term PDF of the linear stochastic model with white noise is Gaussian, described by the mean and the variance, exactly due to the underlying linearity of the dynamics and additive noise. The linear stochastic model is uniquely determined by the three parameters $a$, $f$, and $\sigma$. Three conditions are needed to calibrate these three parameters. Here, these three conditions can be given by the three following statistics: equilibrium mean $\mu$, equilibrium variance $R$, and decorrelation time $\tau$, where $\tau$ is defined as:
\begin{equation}\label{Decorrelation_SI}
    \tau = \int_0^\infty ACF(s)ds,
\end{equation}
with $ACF(s)$ being the autocorrelation function, defined by,
\begin{equation}\label{ACF_SI}
    ACF(s) = \frac{E[(x_t-\mu)(x_{t+s}-\mu)]}{R},
\end{equation}
where the expectation is taken over time. Note that the autocorrelation function \eqref{ACF_SI} is a second-order statistic in time while the equilibrium mean and variance characterize the first- and the second-order long-term statistics. Therefore, the linear stochastic model is expected to capture up to the second-order spatiotemporal statistics of the given time series and perform similarly under such criteria (known as model memory and model fidelity, respectively). The relationship between $(a, f, \sigma)$ and $(\mu,R,\tau)$ is as follows,
\begin{equation}\label{LSM_Calibration_SI}
  \mu = \frac{f}{a},\qquad R = \frac{\sigma^2}{2a},\qquad\mbox {and}\qquad \tau = \frac{1}{a}.
\end{equation}

The procedure of determining the linear stochastic model is as follows:
\begin{itemize}
  \item Step 1. Obtain the training time series. This time series can come from direct observational data or a model, either a direct simulation of a (stochastic) ODE or one of the spectral modes by decomposing a PDE model.
  \item Step 2. Compute the three statistics $(\mu,R,\tau)$ from the time series obtained in Step 1. The mean and the variance are computed using the formulae in \eqref{Moments_SI}. The decorrelation time is calculated from \eqref{Decorrelation_SI}--\eqref{ACF_SI}.
  \item Step 3. Use the relationships in \eqref{LSM_Calibration_SI} to compute the three parameters $(a,f,\sigma)$.
  \item Step 4. Once the parameters in the linear stochastic model are determined, the model can be run for a simulation, which allows the validation of the model's performance.
\end{itemize}

Note that when the given time series is complex-valued, the linear stochastic model will contain one more parameter that represents the phase term describing the oscillations between the real and the imaginary components of the time series. Correspondingly, the cross-correlation function between the real and the imaginary components can be used as the fourth statistic for model calibration that allows us to uniquely determine the four parameters in the complex-valued linear stochastic model (i.e., use the real and imaginary parts of the decorrelation time along with the equilibrium mean and equilibrium variance).

\paragraph{Example.} Consider the following equation for generating a true observed time series
\begin{equation}\label{Cubic_Model_SI}
    \frac{\d x}{\d t} = (f_N +a_Nx+b_Nx^2 -c_Nx^3)+(A_N-B_Nx)\dot{W}_M +\sigma_N \dot{W}_A.
\end{equation}
This is a canonical one-dimensional reduced climate model for low-frequency atmospheric variability with cubic nonlinearities and correlated additive and multiplicative stochastic forcing. Its usages include an application in a regression strategy
using data from a prototype climate model to build one-dimensional stochastic models for
low-frequency patterns such as the North Atlantic Oscillation. Such systematically reduced stochastic models are attractive for sensitivity studies of the climate and other complex high-dimensional systems since they are computationally efficient with comparable prediction skill, and allow for a better understanding of the climate due to their reduced complexity\footnote{Majda, Andrew J., Christian Franzke, and Daan Crommelin. ``Normal forms for reduced stochastic climate models." Proceedings of the National Academy of Sciences 106.10 (2009): 3649-3653.}. The model has several parameters controlling the nonlinearity and state dependent noise. Therefore, the model can exhibit different dynamical and statistical behavior. Four different dynamical regimes are considered here for the true signal. Then the linear stochastic model is calibrated based on the three statistics described above. We aim to see the performance of the linear stochastic model in different regimes. See Table \ref{Table:Cubic_Model_Regimes} for the parameter values in each dynamical regime.

\begin{table}[htb]\centering
\begin{tabular}{|l|c|c|c|c|c|c|c|}
  \hline
    & $a_N$ & $b_N$ & $c_N$ & $f_N$ & $A_N$ & $B_N$ & $\sigma_N$ \\\hline
  Nearly Gaussian regime & -2.2 & 0 & 0 & 2 & 0.1 & 0.1 & 1 \\
  Highly skewed regime & -4 & 2 & 1 & 0.1 & 1 & -1 & 1 \\
  Fat-tailed regime & -3 & -1.5 & 0.5 & 0 & 0.5 & -1 & 1 \\
  Bimodal regime & 4 & 2 & 1 & 0.1 & 1 & -1 & 1 \\
  \hline
\end{tabular}\caption{Dynamical regimes and parameter values of the model \eqref{Cubic_Model_SI}.}\label{Table:Cubic_Model_Regimes}
\end{table}

\begin{figure}[ht!]
	\begin{center}
		\hspace*{-0.5cm}\includegraphics[width=17.5cm]{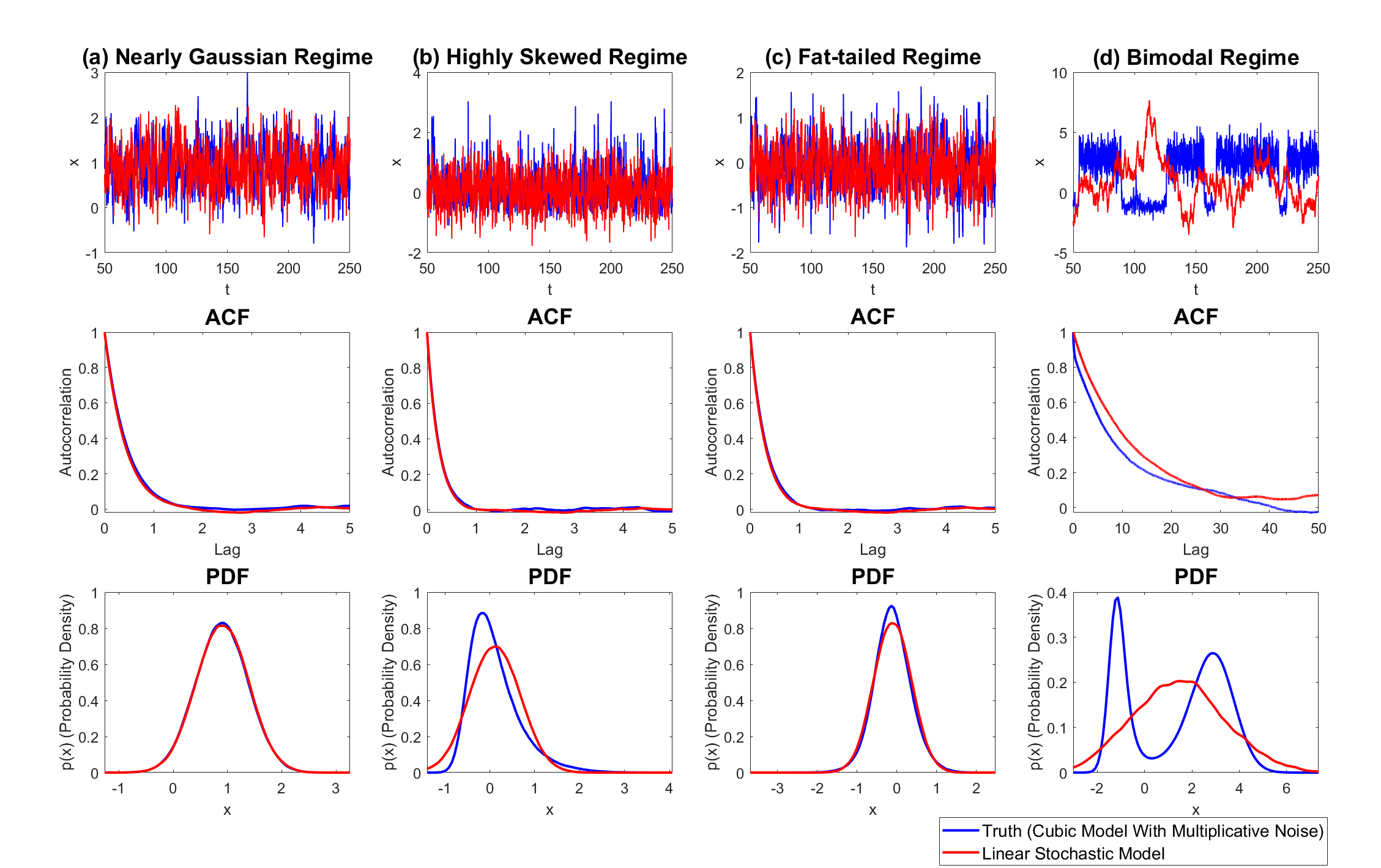}
	\end{center}
	\caption{Comparison of the calibrated linear stochastic model with the nonlinear model \eqref{Cubic_Model_SI}. The four columns show the results in the four dynamical regimes listed in Table \ref{Table:Cubic_Model_Regimes}. The first row compares the time series from the two models. Note that time series of the two models are generated using different random numbers. Therefore, the two time series are not expected to have point-to-point match. The second and the third rows compare the autocorrelation functions and the PDFs, respectively.   }
	\label{SI_LSM}
\end{figure}

Figure \ref{SI_LSM} compares the calibrated linear stochastic model with the nonlinear model \eqref{Cubic_Model_SI}. A long time series with $5000$ time units is utilized for the true signal, avoiding the error from insufficient data. The influence of using short time series on the accuracy of model calibration is a separate UQ topic. When the true signal comes from a nearly Gaussian regime (Column (a)), the linear stochastic model can also perfectly capture the dynamical and statistical behavior. This is not surprising since the first two moments are, by design, the same as the truth. The calibrated linear stochastic model still performs well in highly skewed and fat-tailed regimes (Columns (b)--(c)). The Gaussian PDFs from the linear stochastic model best fit the truth considering only the first two moments. The missing information comes from the extreme events corresponding to the fat tails in the PDFs, which the linear stochastic model cannot capture effectively. The linear stochastic model does not perform well in the bimodal regime (Column (d)) since the truth is too far from the dynamics of the linear stochastic model. A Gaussian distribution cannot approximate the bimodal statistics well. In addition, the state switching in the time series due to the bimodal statistics is missed in the linear stochastic model. In practice, other simple approximate models can reproduce such behavior, such as the linear model with a simple multiplicative noise \cite{chen2023stochastic}.

\begin{comment}
\lstset{language=Matlab,%
    %basicstyle=\color{red},
    breaklines=true,%
    morekeywords={matlab2tikz},
    keywordstyle=\color{blue},%
    morekeywords=[2]{1}, keywordstyle=[2]{\color{black}},
    identifierstyle=\color{black},%
    stringstyle=\color{mylilas},
    commentstyle=\color{mygreen},%
    showstringspaces=false,%without this there will be a symbol in the places where there is a space
    %numbers=left,%
    %numberstyle={\tiny \color{black}},% size of the numbers
    %numbersep=9pt, % this defines how far the numbers are from the text
    emph=[1]{for,end,break},emphstyle=[1]\color{red}, %some words to emphasise
    %emph=[2]{word1,word2}, emphstyle=[2]{style},
}

\section{Matlab Code}
\subsection{Computing the Shannon's entropy}
\lstinputlisting{codes_final/Computing_Entropy.m}
\subsection{Computing the relative entropy}
\lstinputlisting{codes_final/Computing_Relative_Entropy.m}
\end{comment}
\end{document}